\documentclass[aos,MSNbibl,nameyear,noautosecdot,dvips]{arximspdf}
\usepackage{multirow}
\usepackage{graphicx}

\doi{10.1214/12-AOS1041} 
\volume{40}
\issue{5}
\pubyear{2012}
\firstpage{2541}
\lastpage{2571}

\makeatletter

\newtheorem{theorem}{Theorem}

\newtheorem{lemma}{Lemma}

\newtheorem{corollary}{Corollary}

\newcommand\independent{\perp\!\!\!\!\perp}

\def\sign{\operatorname{sign}}

\def\tr{\operatorname{tr}}

\def\aset{\mathcal{A}}

\def\subaset{\mathcal{A}_k}
\def\bolde{\mathbf{e}}
\def\boldf{\mathbf{f}}
\def\bv{\mathbf{v}}
\def\bx{\mathbf{x}}
\def\bw{\mathbf{w}}

\def\br{\mathbf{r}}

\def\bA{\mathbf{A}}
\def\bI{\mathbf{I}}
\def\bH{\mathbf{H}}

\def\hbr{\hat{\mathbf{r}}}
\def\hbR{\hat{\mathbf{R}}}
\def\bM{\mathbf{M}}
\def\bX{\mathbf{X}}

\def\bZ{\mathbf{Z}}
\def\bz{\mathbf{z}}
\def\hbZ{\hat{\mathbf{Z}}}
\def\tbZ{\tilde{\mathbf{Z}}}
\def\bH{\mathbf{H}}
\def\bW{\mathbf{W}}
\def\bmu{\bolds{\mu}}
\def\bzero{\mathbf{0}}

\def\bSigma{\bolds{\Sigma}}
\def\bsigma{\bolds{\sigma}}
\def\hbsigma{\hat{\bolds{\sigma}}}
\def\hbSigma{\hat{\bolds{\Sigma}}}
\def\bTheta{\bolds{\Theta}}
\def\btheta{\bolds{\theta}}
\def\hbtheta{\hat{\bolds{\theta}}}

\def\bbTheta{\breve{\bolds{\Theta}}}
\def\hbTheta{\hat{\bolds{\Theta}}}
\def\tbTheta{\tilde{\bolds{\Theta}}}

\def\balpha{\bolds{\alpha}}
\def\balphaplus{\bolds{\alpha}^+_k}
\def\balphaminus{\bolds{\alpha}^-_k}
\def\talphaplus{\tilde{\bolds{\alpha}}_k^+}
\def\talphaminus{\tilde{\bolds{\alpha}}_k^-}
\def\brbetaplus{\breve{\bolds{\beta}}_k^+}
\def\brbetaminus{\breve{\bolds{\beta}}_k^-}
\def\talpha{\tilde{\bolds{\alpha}}}
\def\talphaplus{\tilde{\bolds{\alpha}}_k^+}
\def\talphaminus{\tilde{\bolds{\alpha}}_k^-}
\def\bbeta{\bolds{\beta}}
\def\brbeta{\breve{\bolds{\beta}}}
\def\hbeta{\hat{\bolds{\beta}}}
\def\tbeta{\tilde{\bolds{\beta}}}
\def\E{\mathrm{E}}
\def\Var{\operatorname{Var}}

\newcommand{\eqref}[1]{(\ref{#1})}
\makeatother

\begin{document}
\begin{frontmatter}

\title{Regularized rank-based estimation of high-dimensional
nonparanormal graphical models\thanksref{T1}}
\runtitle{Rank-based estimation of nonparanormal graphical models}
\thankstext{T1}{Some results in this paper were reported in a
research proposal funded by the Office of Naval Research in November
2010. Supported in part by NSF Grant DMS-08-46068 and a grant from ONR.}

\begin{aug}
\author{\fnms{Lingzhou} \snm{Xue}\ead[label=e1]{lzxue@stat.umn.edu}}
\and
\author{\fnms{Hui} \snm{Zou}\corref{}\ead[label=e2]{zouxx019@umn.edu}}
\runauthor{L. Xue and H. Zou}
\affiliation{University of Minnesota}
\address{School of Statistics\\
University of Minnesota\\
Minneapolis, Minnesota 55455\\
USA\\
\printead{e1}\\
\phantom{E-mail:\ }\printead*{e2}} 
\end{aug}

\received{\smonth{2} \syear{2012}}
\revised{\smonth{8} \syear{2012}}

%
\begin{abstract}
A sparse precision matrix can be directly translated into a sparse
Gaussian graphical model under the assumption that the data follow
a joint normal distribution. This neat property makes high-dimen\-sional
precision matrix estimation very appealing in many applications.
However, in practice we often face nonnormal data, and variable
transformation is often used to achieve normality. In this paper we
consider the nonparanormal model that assumes that the variables follow
a joint normal distribution after a set of unknown monotone
transformations. The nonparanormal model is much more flexible than the
normal model while retaining the good interpretability of the latter in
that each zero entry in the sparse precision matrix of the
nonparanormal model corresponds to a pair of conditionally independent
variables. In this paper we show that the nonparanormal graphical model
can be efficiently estimated by using a rank-based estimation scheme
which does not require estimating these unknown transformation
functions. In particular, we study the rank-based graphical lasso, the
rank-based neighborhood Dantzig selector and the rank-based CLIME. We
establish their theoretical properties in the setting where the
dimension is nearly exponentially large relative to the sample size.
It is shown that the proposed rank-based estimators work as well as
their oracle counterparts defined with the oracle data. Furthermore,
the theory motivates us to consider the adaptive version of the
rank-based neighborhood Dantzig selector and the rank-based CLIME that
are shown to enjoy graphical model selection consistency without
assuming the irrepresentable condition for the oracle and rank-based
graphical lasso. Simulated and real data are used to demonstrate the
finite performance of the rank-based estimators.
\end{abstract}

%

\end{frontmatter}

\section{\texorpdfstring{Introduction.}{Introduction}}\label{sec1}

Estimating covariance or precision matrices is of fundamental
importance in multivariate statistical methodologies and applications.
In particular, when data follow a joint normal distribution, $
\bX=(X_1,\ldots,X_p)\sim N_{p}(\bmu,\bSigma),
$ the precision matrix $\bTheta=\bSigma^{-1}$ can be directly
translated into a Gaussian graphical model. The Gaussian graphical
model serves as a noncausal structured approach to explore the complex
systems consisting of Gaussian random variables, and it finds many
interesting applications in areas such as gene expression genomics and
macroeconomics determinants study
[\citet{nfriedman2004,wille2004,dobra2009}]. The precision matrix plays
a critical role in the Gaussian graphical models
because the zero entries in $\bTheta= (\theta_{ij} )_{p\times
p}$ precisely capture the desired conditional independencies, that is,
$
\theta_{ij}=0
$
if and only if
$
X_i \independent X_j | \bX\setminus\{X_i,X_j\}$ [\citet
{lauritzen1996,edwards2000}].

The sparsity pursuit in precision matrices was initially considered by
\citet{dempster1972} as the covariance selection problem.
Multiple testing methods have been employed for network exploration in
the Gaussian graphical models [\citet{drton2004}]. With rapid advances
of the high-throughput technology (e.g., microarray, functional MRI),
estimation of a sparse graphical model has become increasingly
important in the high-dimensional setting. Some well-developed
penalization techniques have been used for estimating sparse Gaussian
graphical models. In a highly-cited paper, \citet{meinshausen2006}
proposed the neighborhood selection scheme which tries to discover the
smallest index set $ne_{\alpha}$ for each variable $X_{\alpha}$
satisfying $X_{\alpha} \independent\bX\setminus\{X_{\alpha},\bX_{ne_{\alpha}}\}| \bX_{ne_{\alpha}}$. \citet{meinshausen2006} further
proposed to use the lasso [\citet{tibshirani1996}] to fit each
neighborhood regression model. Afterwards, one can summarize the zero
patterns by aggregation via union or intersection. \citet{yuan2010}
considered the Dantzig selector [\citet{candes2007}] as an alternative
to the lasso penalized least squares in the neighborhood selection
scheme. \citet{peng2009} proposed the joint neighborhood lasso
selection. Penalized likelihood methods have been studied for Gaussian
graphical modeling [\citet{yuan2007}]. \citet{friedman2008} developed a
fast blockwise coordinate descent algorithm [\citet{banerjee2008}]
called graphical lasso for efficiently solving the lasso penalized
Gaussian graphical model. 
Rate of convergence under the Frobenius norm was established by \citet
{rothman2008}. \citet{ravikumar2008} obtained the convergence rate under
the elementwise $\ell_{\infty}$ norm and the spectral norm. \citet
{fan2009a} studied the nonconvex penalized Gaussian graphical model
where a nonconvex penalty such as SCAD [\citet{fan2001}] is used to
replace the lasso penalty in order to overcome the bias issue of the
lasso penalization. \citet{zhou2010} proposed a hybrid method for
estimating sparse Gaussian graphical models: they first infer a sparse
Gaussian graphical model structure via thresholding neighborhood
selection and then estimate the precision matrix of the submodel by
maximum likelihood. \citet{clime11} recently proposed a constrained $\ell_1$
minimization estimator called CLIME for estimating sparse precision
matrices and established its convergence rates under the elementwise
$\ell_{\infty}$ norm and Frobenius norm.

%
\begin{table}
\caption{Testing for normality of the gene expression measurements in
the Arabidposis thaliana
data. This table illustrates the number out of 39 genes rejecting
the null hypothesis
of normality at the significance level of $0.05$}\label{realdatanormalitytest}
\begin{tabular*}{\textwidth}{@{\extracolsep{\fill}}lcccc@{}}
\hline
& \textbf{Critical value}& \textbf{Cramer--von Mises} & \textbf{Lilliefors} & \multicolumn{1}{c@{}}{\textbf{Shapiro--Francia}} \\
\hline
{Raw data} & 0.05 &30 & 30 & 35 \\
& 0.05$/$39 &24 & 26 & 28 \\
{Log data} & 0.05 &29 & 24 & 33 \\
& 0.05$/$39 &14 & 12 & 16 \\
\hline
\end{tabular*}
\end{table}

\begin{figure}

\includegraphics{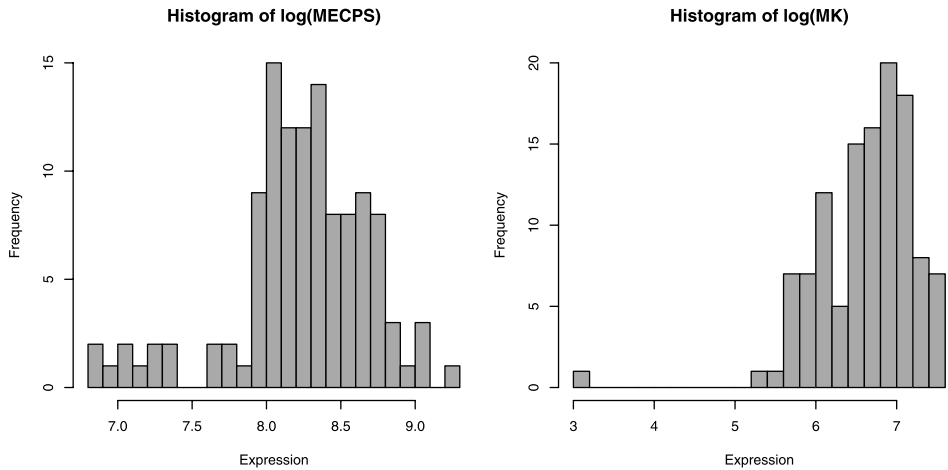}

\caption{Illustration of the nonnormality after the log-transformation
preprocessing.}\label{realdatanormalityhist}
\end{figure}

Although the normality assumption can be relaxed if we only focus on
estimating a precision matrix, it plays an essential role in making the
neat connection between a sparse precision matrix and a sparse Gaussian
graphical model. Without normality, we ought to be very cautious when
translating a good sparse precision matrix estimator into an
interpretable sparse Gaussian graphical model. However, the normality
assumption often fails in reality. For example, the observed data are
often skewed or have heavy tails. To illustrate the issue of
nonnormality in real applications, let us consider the gene expression
data to construct isoprenoid genetic regulatory network in Arabidposis
thaliana [\citet{wille2004}], including 16 genes from the mevalonate
(MVA) pathway in the cytosolic, 18 genes from the plastidial (MEP)
pathway in the chloroplast and 5 encode proteins in the mitochondrial.
This dataset contains gene expression measurements of 39 genes assayed
on $n = 118$ Affymetrix GeneChip microarrays. 
This dataset was analyzed by \citet{wille2004}, \citet{lihz2006} and \citet
{drton2007} in the context of Gaussian graphical modeling after taking
the log-transformation of the data. However, the normality assumption
is still inappropriate even after the log-transformation. To show this,
we conduct the normality test at the significance level of $0.05$ as in
Table~\ref{realdatanormalitytest}. It is clear that at most 9 out of
39 genes would pass any of three normality tests. Even after
log-transformation, at least $60\%$ genes reject the null hypothesis of
normality. With Bonferroni correction there are still over $30\%$ genes
that fail to pass any normality test. Figure \ref
{realdatanormalityhist} plots the histograms of two key isoprenoid
genes MECPS in the MEP pathway and MK in the MVA pathway after the
log-transformation, clearly showing the nonnormality of the data after
the log-transformation.

Using transformation to achieve normality is a classical idea in
statistical modeling. The celebrated Box--Cox transformation is widely
used in regression analysis. However, any parametric modeling of the
transformation suffers from model mis-specification which could lead to
misleading inference. In this paper we take a nonparametric
transformation strategy to handle the nonnormality issue. Let $F(\cdot
)$ be the CDF of a continuous random variable $X$ and $\Phi^{-1}(\cdot
)$ be the inverse of the CDF of $N(0,1)$. Consider the transformation
from $X$ to $Z$ by $Z=\Phi^{-1}(F(X))$. Then it is easy to see that $Z$
is standard normal regardless of $F$. Motivated by this simple fact, we
consider modeling the data by the following nonparanormal model:\vspace*{8pt}

\textit{The nonparanormal model}: $\bX=(X_1,\ldots,X_p)$
follows a $p$-dimensional nonparanormal distribution if there exists a
vector of unknown univariate monotone increasing transformations,
denoted by $\boldf=(f_1,\ldots,f_p)$, such that the transformed random
vector follows a multivariate normal distribution with mean 0 and
covariance $\bSigma$,
\begin{equation}
\label{model} \boldf(\bX)= \bigl(f_1(X_1),
\ldots,f_p(X_p) \bigr) \sim N_{p}(0,\bSigma),
\end{equation}
where without loss of generality the diagonals of $\bSigma$ are equal
to 1.\vspace*{8pt}

Note that model (\ref{model}) implies that $f_j(X_j)$ is a standard
normal random variable. Thus,
$f_j$ must be $\Phi^{-1}\circ F_j$ where $F_j$ is the CDF of $X_j$.
The marginal normality is always achieved by transformations, so model
(\ref{model}) basically assumes that those marginally
normal-transformed variables are jointly normal as well.
We follow \citet{lafferty2009} to call model (\ref{model}) the
nonparanormal model, but model (\ref{model}) is in fact a
semiparametric Gaussian copula model. { The semiparametric Gaussian
copula model is a nice combination of flexibility and interpretability,
and it has generated a lot of interests in statistics, econometrics and
finance; see \citet{song2000}, \citet{chen2006,chen2006jasa} and
references therein. Let $ \bZ=(Z_1,\ldots,Z_p)= (f_1(X_1),\ldots
,f_p(X_p) )$. By the joint normality assumption of $\bZ$, we know
that $\theta_{ij}=0$ if and only if $Z_i \independent Z_j |  \bZ
\setminus\{Z_i,Z_j\}.$ Interestingly, we have that
\[
Z_i \independent Z_j |  \bZ\setminus
\{Z_i,Z_j\} \quad\Longleftrightarrow\quad X_i
\independent X_j |  \bX\setminus\{X_i,X_j
\}.
\]
Therefore, a sparse $\bTheta$ can be directly translated into a sparse
graphical model for presenting the original variables.

In this work we primarily focus on estimating $\bTheta$ which is then
used to construct a nonparanormal graphical model. As for the
nonparametric transformation function, by the expression $f_j=\Phi^{-1}\circ F_j$,
we have a natural estimator for the transformation
function of the $j$th variable as $\hat f_j=\Phi^{-1}\circ\hat
F^{+}_j$ where $\hat F^{+}_j$ is a Winsorized empirical CDF of the
$j$th variables. Note that the Winsorization is used to avoid infinity
value and to achieve better bias-variance tradeoff; see \citet
{lafferty2009} for detailed discussion.
In this paper we show that we can directly estimate $\bTheta$ without
estimating these nonparametric transformation functions at all. This
statement seems to be a bit surprising because a natural estimation
scheme is a two-stage procedure: first estimate $f_j$ and then apply a
well-developed sparse Gaussian graphical model estimation method to the
transformed data $\hat{\bz}_i=\mathbf{\hat f}(\bx_i), 1 \le i \le
n$. \citet{lafferty2009} have actually studied this ``plug-in''
estimation approach. They proposed a Winsorized estimator of the
nonparametric transformation function and used the graphical lasso in
the second stage. They established convergence rate of the ``plug-in''
estimator when $p$ is restricted to a polynomial order of $n$. However,
\citet{lafferty2009} did not get a satisfactory rate of convergence for
the ``plug-in'' approach, because the rate of convergence can be
established for the Gaussian graphical model even when $p$ grows with
$n$ almost exponentially fast [\citet{ravikumar2008}]. As noted in \citet
{lafferty2009}, it is very challenging, if not impossible, to push the
theory of the ``plug-in'' approach to handle exponentially large dimensions.
One might ask if using a better estimator for the transformation
functions could improve the rate of convergence such that $p$ could be
allowed to be nearly exponentially large relative to $n$. This is a
legitimate direction for research. We do not pursue this direction in
this work. Instead, we show that we could use a rank-based estimation
approach to achieve the exact same goal without estimating these
transformation functions at all.

Our estimator is constructed in two steps. First, we propose using the
adjusted Spearman's rank correlation to get a nonparametric sample
estimate of $\bSigma$. As the second step, we compute a sparse
estimator $\bTheta$ from the rank-based sample estimate of $\bSigma$.
For that purpose, we consider several regularized rank estimators,
including the \textit{rank-based graphical lasso}, the \textit{rank-based
neighborhood Dantzig selector} and the \textit{rank-based CLIME}. The
complete methodological details are presented in Section~\ref{sec2}.
In Section~\ref{sec3} we establish theoretical properties of the proposed rank-based
estimators, regarding both precision matrix estimation and graphical
model selection. In particular, we are motivated by the theory to
consider the adaptive version of the rank-based neighborhood Dantzig
selector and the rank-based CLIME, which can select the true support
set with an overwhelming probability without assuming a stringent
irrepresentable condition required for the oracle and rank-based
graphical lasso. Section~\ref{sec4} contains numerical results and Section~\ref{sec5}
has some concluding remarks. Technical proofs are presented in an \hyperref[app]{Appendix}.

A referee informed us in his/her review report that \citet{liu2012}
also independently used the rank-based correlation in the context of
nonparametric Gaussian graphical model estimation. A major focus in Liu
et al. (\citeyear{liu2012}) is the numerical demonstration of the robustness property
of the rank-based methods using both Spearman's rho and Kendall's tau
when data are contaminated. In the present paper we provide a
systematic analysis of the rank-based estimators, and our theoretical
analysis further leads to the rank-based adaptive Dantizg selector and
the rank-based adaptive CLIME in order to achieve improved sparsity
recovery properties. Our theoretical analysis of the rank-based
adaptive Dantizg selector is of independent interest. Although the
theory is established for the rank-based estimators using Spearman's
rho, the same analysis can be easily adopted to prove the theoretical
properties of the rank-based estimators using Kendall's tau rank correlation.

\section{\texorpdfstring{Methodology.}{Methodology}}\label{sec2}

We first introduce some necessary notation. For a matrix $\bA
=(a_{ij})$, we define its entry-wise $\ell_1$ norm as $\|\bA\|_1=\sum_{(i,j)}|a_{ij}|$,
and its entry-wise $\ell_{\infty}$ norm as $\|\bA\|_{\max}=\max_{(i,j)}|a_{ij}|$.
For a vector $\bv=(v_1,\ldots,v_l)$,
we define its $\ell_1$ norm as $\|\bv\|_{\ell_1}=\sum_j|v_j|$ and its
$\ell_{\infty}$ norm as $\|\bv\|_{\ell_\infty}=\max_j|v_j|$. To
simplify notation, define $\bM_{A,B}$ as the sub-matrix of $\bM$ with
row indexes $A$ and column indexes $B$, and define $\bv_A$ as the
sub-vector of $\bv$ with indexes $A$. Let $(k)$ be the index set $\{
1,\ldots,k-1,k+1,\ldots,p\}$. Denote by { $\bSigma_{(k)}=\bSigma_{(k),(k)}$}
the sub-matrix of $\bSigma$ with both $k$th row and column
removed, and denote by $\bsigma_{(k)}=\bSigma_{(k),k}$ the vector
including all the covariances associated with the $k$th variable. In
the same fashion, we can also define $\bTheta_{(k)}$, $\btheta_{(k)}$,
and so on.

\subsection{\texorpdfstring{The ``oracle'' procedures.}{The ``oracle'' procedures}}\label{sec2.1}
Suppose an oracle knows the underlying transformation vector; then the
oracle could easily recover ``oracle data'' by applying these true
transformations, that is, $\bz_i=\boldf(\bx_i), 1 \le i \le n$.
Before presenting our rank-based estimators, it is helpful to revisit
the ``oracle'' procedures that are defined based on the ``oracle data.''

\begin{itemize}
\item\textit{The oracle graphical lasso.} Let $\hbSigma^o$ be the
sample covariance matrix for the ``oracle'' data, and then the
``oracle'' log-profile-likelihood becomes $\log\det(\bTheta)-\tr
(\hbSigma^o\bTheta)$. The ``oracle'' graphical lasso solves the
following $\ell_1$ penalized likelihood problem:
\begin{equation}
\label{sec21eq1} \min_{\bTheta\succ0}-\log\det(\bTheta)+\tr \bigl(
\hbSigma^o\bTheta \bigr)+\lambda \sum_{i\neq j}
|\theta_{ij}|.
\end{equation}

\item\textit{The oracle neighborhood lasso selection.} Under the
nonparanormal model, for each $k=1,\ldots,p$, the ``oracle'' variable
$Z_k$ given $\bZ_{(k)}$ is normally distributed as $
N (\bZ_{(k)}^T\bSigma_{(k)}^{-1}\bsigma_{(k)}, 1-\bsigma_{(k)}^T\bSigma_{(k)}^{-1}\sigma_{(k)} ),$ which can be written as $
Z_k=\bZ_{(k)}^T\bbeta_{k}+\varepsilon_k
$ with $ \bbeta_{k}=\bSigma_{(k)}^{-1}\bsigma_{(k)}$ and $
\varepsilon_k\sim N(0, 1-\bsigma_{(k)}^T\bSigma_{(k)}^{-1}\bsigma_{(k)}).
$\vspace*{1pt} Notice that $\bbeta_k$ and $\varepsilon_k$ are closely related to
the precision matrix $\bTheta$, that is, $
\theta_{kk}=1/\Var(\varepsilon_k)
$ and $
\btheta_{(k)}=-\bbeta_{k}/\Var(\varepsilon_k).
$
Thus for the $k$th variable, $\btheta_{(k)}$ and $\bbeta_k$ share the
same sparsity pattern. Following \citet{meinshausen2006}, the oracle
neighborhood lasso selection obtains the solution $\hbeta^o_k$ from the
following lasso penalized least squares problem:
\begin{equation}
\label{sec21eq2} \min_{\bbeta\in\mathbb{R}^{p-1}} \frac{1}n\sum
_{i=1}^n \bigl(z_{ik}-
\bz_{i(k)}^T\bbeta \bigr)^2+\lambda\|\bbeta
\|_{\ell_1},
\end{equation}
and then the sparsity pattern of $\bTheta$ can be estimated by
aggregating the neighborhood support set of $\hbeta^o_k=(\hat\beta^o_{jk})_{j\neq k}$ ($\widehat{ne}_{k}=\{j\dvtx  \hat\beta^o_{jk}\neq0\}$) via
intersection or union.

We notice the fact that
\[
\frac{1}n\sum_{i=1}^n{
\bigl(}z_{ik}-\bz_{i(k)}^T\bbeta{
\bigr)}^2 = \bbeta^T\hbSigma^o_{(k)}
\bbeta-2\bbeta^T\hbsigma^o_{(k)}+\hat
\sigma^o_{kk}.
\]
Then (\ref{sec21eq2}) can be written in the following equivalent form:
\begin{equation}
\label{sec21eq2b} \min_{\bbeta\in\mathbb{R}^{p-1}} \bbeta^T\hbSigma^o_{(k)}
\bbeta-2\bbeta^T\hbsigma^o_{(k)}+\lambda\|
\bbeta \|_{\ell_1}.
\end{equation}

\item\textit{The oracle neighborhood Dantzig selector.} Following \citet
{yuan2010} the lasso least squares in (\ref{sec21eq2}) can be replaced
with the Dantzig selector
\begin{equation}
\label{sec21eq3} \min_{\bbeta\in\mathbb{R}^{p-1}}\|\bbeta\|_{\ell_1} \qquad\mbox{subject
to } \Biggl\|\frac{1}n\sum_{i=1}^n
\bz_{i(k)} \bigl(\bz_{i(k)}^T\bbeta -z_{ik}
\bigr)\Biggr\|_{\ell_\infty}\le\lambda.
\end{equation}
Then the sparsity pattern of $\bTheta$ can be similarly estimated by
aggregating via intersection or union. Furthermore, we notice that
\[
\frac{1}n\sum_{i=1}^n
\bz_{i(k)} \bigl(\bz_{i(k)}^T\bbeta-z_{ik}
\bigr) = \hbSigma^o_{(k)}\bbeta-\hbsigma^o_{(k)}.
\]
Then (\ref{sec21eq3}) can be written in the following equivalent form:
\begin{equation}
\label{sec21eq3b} \min_{\bbeta\in\mathbb{R}^{p-1}}\|\bbeta\|_{\ell_1} \qquad\mbox{subject
to } \bigl\|\hbSigma^o_{(k)}\bbeta-\hbsigma^o_{(k)}
\bigr\|_{\ell_\infty}\le\lambda.
\end{equation}

\item\textit{The oracle CLIME.} Following \citet{clime11} we can
estimate precision matrices by solving a constrained $\ell_1$
minimization problem,
\begin{equation}
\label{sec21eq4a} \mathop{\arg\min}_{\bTheta} \|\bTheta\|_1\qquad
\mbox{subject to } \bigl\|\hbSigma^o\bTheta-\bI\bigr\|_{\max}\le
\lambda.
\end{equation}
\citet{clime11} compared the CLIME with the graphical lasso, and showed
that the CLIME enjoys nice theoretical properties without assuming the
irrepresentable condition of Ravikumar et al. (\citeyear{ravikumar2008}) for the graphical lasso.
\end{itemize}

\subsection{\texorpdfstring{The proposed rank-based estimators.}{The proposed rank-based estimators}}\label{sec2.2}
The existing theoretical results in the literature can be directly
applied to these oracle estimators.
However, the ``oracle data'' $\bz_1,\bz_2,\ldots,\bz_n$ are unavailable
and thus the above-mentioned ``oracle'' procedures are not genuine estimators.
Naturally we wish to construct a genuine estimator that can mimic the
oracle estimator. To this end, we can derive an alternative estimator
of $\bSigma$ based on the actual data $\bx_1,\bx_2,\ldots,\bx_n$ and
then feed this genuine covariance estimator to the graphical lasso, the
neighborhood selection or CLIME. To implement this natural idea, we
propose a rank-based estimation scheme. Note that $\bSigma$ can be
viewed as the correlation matrix as well, that is, $\sigma_{ij}=\operatorname
{corr}(\bz_i,\bz_j)$.
Let ($x_{1i},x_{2i},\ldots,x_{ni}$) be the observed values of variable
$X_i$. We convert them to ranks denoted by
$\br_i=(r_{1i},r_{2i},\ldots,r_{ni})$.
Spearman's rank correlation $\hat r_{ij}$ is defined as Pearson's
correlation between $\br_i$ and $\br_j$.
Spearman's rank correlation is a nonparametric measure of dependence
between two variables. It is important to note
that $\br_i$ are the ranks of the ``oracle'' data.
Therefore, $\hat r_{ij}$ is also identical to the Spearman's rank
correlation between the ``oracle'' variables $Z_i,Z_j$.
In other words, in the framework of rank-based estimation, we can treat
the observed data as the ``oracle'' data and avoid estimating $p$
nonparametric transformation functions.
We make a note here that one may consider other rank correlation
measures such as Kendall's tau correlation. To fix the idea we use
Spearman's rank correlation throughout this paper.

The nonparanormal model implies that $(Z_i,Z_j)$ follows a bivariate
normal distribution with correlation parameter $\sigma_{ij}$.
Then a classical result due to \citet{kendall1948} [see also \citet
{kruskal1958}] shows that
\begin{equation}
\label{bias} \lim_{n\rightarrow+\infty} \E(\hat r_{ij})=\frac{6}{\pi}
\arcsin \biggl(\frac
{1}2\sigma_{ij} \biggr),
\end{equation}
which indicates that $\hat r_{ij}$ is a biased estimator of $\sigma_{ij}$. To correct the bias,
\citet{kendall1948} suggested using the adjusted Spearman's rank correlation
\begin{equation}
\label{adjust} \hat r^s_{ij}=2\sin \biggl(
\frac{\pi}6 \hat r_{ij} \biggr).
\end{equation}
Combining (\ref{bias}) and (\ref{adjust}) we see that $\hat r^s_{ij}$
is an asymptotically unbiased estimator of~$\sigma_{ij}$.
Naturally we define the rank-based sample estimate of $\bSigma$ as follows:
\[
\label{ranksample} \hbR^s= \bigl(\hat r^s_{ij}
\bigr)_{1\le i,j\le p}.
\]

In Section~\ref{sec3} we show $\hbR^s$ is a good estimator of $\bSigma$. Then we
naturally come up with the following rank-based estimators of $\bTheta$
by using the graphical lasso, the neighborhood Dantzig selector and CLIME:
\begin{itemize}
\item\textit{The rank-based graphical lasso}:
\begin{equation}
\label{sec22eq1} \hbTheta_g^s=\mathop{\arg\min}_{\bTheta\succ0}-
\log\det(\bTheta)+\tr \bigl(\hbR^s\bTheta \bigr)+\lambda\sum
_{i\neq j} |\theta_{ij}|.
\end{equation}

\item\textit{The rank-based neighborhood Dantzig selector:} A
rank-based estimate of $\bbeta_k$ can be solved by
\begin{equation}
\label{sec22eq2} \hbeta_{k}^{s.nd}=\mathop{\arg\min}_{\bbeta\in\mathbb{R}^{p-1}}\|
\bbeta\|_{\ell
_1}\qquad \mbox{subject to } \bigl\|\hbR^s_{(k)}
\bbeta-\hbr^s_{(k)}\bigr\|_{\ell_\infty}\le\lambda.
\end{equation}
The support of $\bTheta$ can be estimated from the support of $\hbeta_1^{s.nd},\ldots,\hbeta_p^{s.nd}$ via aggregation by union or
intersection. We can also construct the rank-based precision matrix
estimator $\hbTheta_{nd}^s=(\hat\theta^{s.nd}_{ij})_{1\le i,j\le p}$ with
\[
\hat\theta^{s.nd}_{kk}= \bigl( \bigl(\hbeta_k^{s.nd}
\bigr)^T\hbR_{(k)}^s\hbeta_k^{s.nd}
-2 \bigl(\hbeta_k^{s.nd} \bigr)^T
\hbr^{s.nd}_{(k)}+1 \bigr)^{-1}
\]
and
\[
\hbtheta^{s.nd}_{(k)}=-\hat\theta^{s.nd}_{kk}
\hbeta_k^{s.nd}
\]
($k=1,\ldots,p$). We can symmetrize $\hbTheta_{nd}^s$ by solving the
following optimization problem [\citet{yuan2010}]:
\[
\bbTheta_{nd}^s=\mathop{\arg\min}_{\bTheta: \bTheta=\bTheta'}\bigl\|\bTheta-
\hbTheta_{nd}^s\bigr\|_{\ell_1}.
\]
Theoretical analysis of the rank-based neighborhood Dantzig selector in
Section~\ref{sec3.2} motivated us to consider using the adaptive Dantzig
selector in the rank-based neighborhood estimation in order to achieve
better graphical model selection performance. See Section~\ref{sec3.2} for more details.

\item\textit{The rank-based CLIME}:
\begin{equation}
\label{sec22eq3} \hbTheta_c^s=\mathop{\arg\min}_{\bTheta}\|
\bTheta\|_1 \qquad\mbox{subject to } \bigl\|\hbR^s\bTheta-\bI
\bigr\|_{\max}\le\lambda.
\end{equation}

Let $\bolde_k$'s be the natural basis in $\mathbb{R}^p$. By Lemma~\ref{concentration} in
\citet{clime11} the above optimization problem can be further decomposed
into $p$ subproblems of vector minimization,
\begin{equation}
\label{sec21eq4b} \hbtheta^{s.c}_k=\mathop{\arg\min}_{\btheta\in\mathbb{R}^p}
\|{\btheta}\|_{\ell
_1} \qquad\mbox{subject to } \bigl\|\hbR^s\btheta-
\bolde_k\bigr\|_{\ell_\infty
}\le\lambda,
\end{equation}
for $k=1,\ldots,p$. Then $\hbTheta_c^s$ is exactly equivalent to
$(\hbtheta^{s.c}_1,\ldots,\hbtheta^{s.c}_p)$.
Note that $\hbTheta_c^s$ could be asymmetric. Following \citet{clime11}
we consider
\[
\bbTheta_c^s= \bigl(\breve\theta^{s.c}_{ij}
\bigr)_{1\le i,j\le p}
\]
with $
\breve\theta^{s.c}_{ij}=\hat\theta^{s.c}_{ij}I_{\{|\hat\theta
^{s.c}_{ij}|\le|\hat\theta^{s.c}_{ji}|\}}
+\hat\theta^{s.c}_{ji}I_{\{|\hat\theta^{s.c}_{ij}|>|\hat\theta
^{s.c}_{ji}|\}}.$
In the original paper \citet{clime11} proposed to use hard thresholding
for graphical model selection. Borrowing the basic idea from \citet
{zou2006}, we propose an adaptive version of the rank-based CLIME in
order to achieve better graphical model selection. See Section~\ref{sec3.3} for
more details.
\end{itemize}

\subsection{\texorpdfstring{Rank-based neighborhood lasso?}{Rank-based neighborhood lasso?}}\label{sec2.3}

One might consider the rank-based neighborhood lasso defined as follows:
\begin{equation}
\label{sec22eq4} \min_{\bbeta\in\mathbb{R}^{p-1}} \bbeta^T\hbR^s_{(k)}
\bbeta-2\bbeta^T\hbr^s_{(k)}+\lambda\|\bbeta
\|_{\ell_1}.
\end{equation}
However, there is a technical problem for the above definition. The
Spearman's rank correlation matrix $\hbR$ is always positive
semidefinite, but the adjusted correlation
matrix $\hbR^s$ could become indefinite. To our best knowledge, \citet
{devlin1975} were the first to point out the indefinite issue of the
estimated rank correlation matrix. Here we also use a toy example to
illustrate this point. Consider the $3\times3$ correlation matrix
\[
\bA=\pmatrix{1 & 0.7 & 0\vspace*{2pt}
\cr
0.7 & 1 & 0.7\vspace*{2pt}
\cr
0 & 0.7
& 1}.
\]
Note that $\bA$ is positive-definite with eigenvalues $1.99$, $1.00$
and $0.01$, but $2\sin(\frac{\pi}6\bA)$ becomes indefinite with
eigenvalues $2.01$, $1.00$ and $-0.01$.
The negative eigenvalues will make (\ref{sec22eq4}) an ill-defined
optimization problem.
Fortunately, the positive definite issue does not cause any problem for
the graphical lasso, Dantzig selector and CLIME. Notice that the
diagonal elements of $\hbR^s$ are obviously strictly positive, and thus
Lemma 3 in \citet{ravikumar2008} suggests that the rank-based graphical
lasso always has a unique positive definite solution { for any
regularization parameter $\lambda>0$}. The rank-based neighborhood
Dantzig selector and the rank-based CLIME are still well defined, even
when $\hbR^s_{(k)}$ becomes indefinite, and the according optimization
algorithms also tolerate the indefiniteness of $\hbR^s_{(k)}$. One
might consider a positive definite correction of $\hbR^s$ for
implementing the neighborhood lasso estimator. However, the resulting
estimator shall behave similarly to the rank-based neighborhood Dantzig
selector because the lasso penalized least squares and Dantzig
selector, in general, work very similarly [\citet{dasso1,dasso2009}].


\section{\texorpdfstring{Theoretical properties.}{Theoretical properties}}\label{sec3}

For a vector $\bv=(v_1,\ldots,v_l)$, let $\|\bv\|_{\min}$ denote the
minimum absolute value, that is, $\|\bv\|_{\min}=\min_j|v_j|$. For a
matrix $\bA=(a_{ij})_{k\times l}$, we define the following matrix
norms: the matrix $\ell_1$ norm $\|\bA\|_{\ell_1}=\max_j\sum_i|a_{ij}|$, the matrix
$\ell_\infty$ norm $\|\bA\|_{\ell_\infty}=\max_i\sum_j|a_{ij}|$ and the Frobenius
norm $\|\bA\|_F=({\sum_{(i,j)}a^2_{ij}})^{1/2}$. For any symmetric matrix, its matrix $\ell_1$
norm coincides its matrix $\ell_\infty$ norm. Denote by
$\lambda_{\min}(\bA)$ and $\lambda_{\max}(\bA)$ the smallest and
largest eigenvalues of $\bA$, respectively. Define $\bSigma^*$ as the
true covariance matrix, and let $\bTheta^*$ be its inverse. Let $\aset
$ be the true support set of the off-diagonal elements in~$\bTheta^*$.
Let $d = \max_j \sum_{i\neq j} I_{\{\theta^*_{ij}\neq0\}}$ be the
maximal degree over the underlying graph corresponding to $\bTheta^*$,
and let $s=\sum_{(i,j)\dvtx i\neq j}I_{\{\theta^*_{ij}\neq0\}}$ be the
total degree over the whole graph.

In this section we establish theoretical properties for the proposed
rank-based estimators. The main conclusion drawn from these theoretical
results is that the rank-based graphical lasso, neighborhood Dantzig
selector and CLIME work as well as their oracle counterparts in terms
of the rates of convergence. We first provide useful concentration
bounds concerning the accuracy of the rank-based sample correlation matrix.

\begin{lemma} \label{concentration}
Fix any $0<\varepsilon<1$, and let $n\ge\frac{12\pi}{\varepsilon}$.
Then there exists some absolute constant $c_0>0$, and we have the
following concentration bounds:
\begin{eqnarray*}
\Pr \bigl(\bigl|\hat r^s_{ij}-\sigma_{ij}\bigr|>
\varepsilon \bigr)&\le& 2\exp \bigl(-c_0n\varepsilon^2
\bigr),
\\
\Pr \bigl(\bigl\Vert\hbR^s-\bSigma\bigr\Vert_{\rm{max} }>\varepsilon
\bigr)& \le& p^2\exp \bigl(-c_0n\varepsilon^2
\bigr).
\end{eqnarray*}
\end{lemma}

Lemma~\ref{concentration} is a key ingredient of our theoretical analysis. It basically
shows that the rank-based sample estimator of $\bSigma$ works as well
as the usual sample covariance estimator of $\bSigma$ based on the
``oracle data.''

\subsection{\texorpdfstring{Rank-based graphical lasso.}{Rank-based graphical lasso}}\label{sec3.1}

Denote by $\psi_{\min}=\min_{(i,j)\in\aset}|\theta^*_{ij}|$ the minimal
entry of $\bTheta^*$ in the absolute scale. Define $K_{\bSigma^*}=\|
\bSigma^*_{\aset\aset}\|_{\ell_{\infty}}$ and $K_{\bH^*}=\|(\bH^*_{\aset
\aset})^{-1}\|_{\ell_{\infty}}$. Define $\bH^*$ as the Kronecker
product $\bSigma^*\otimes\bSigma^*$.

\begin{theorem} \label{nonasymrglasso}
Assume $\|\bH^*_{\aset^c \aset}(\bH^*_{\aset\aset})^{-1}\|_{\ell_{\infty
}}<1-\kappa$ for $\kappa\in(0,1)$.
\begin{longlist}[(a)]
\item[(a)] Element-wise maximal bound: if $\lambda$ is chosen such that
\[
\lambda < \frac{1}{6(1+{\kappa}/4)K_{\bSigma^*}K_{\bH^*}
\max\{1,(1+{4}/{\kappa})K_{\bSigma^*}^2K_{\bH^*}\}}\cdot\frac{1}d,
\]
with probability at least $1-p^2\exp(-\frac{\kappa^2}{16}c_0n\lambda^2)$, the rank-based \textit{graphical lasso} estimator $\hbTheta_g^s$
satisfies that $\hat\theta^{s.g}_{ij}=0$ for any $(i,j)\in\aset^c$ and
\[
\bigl\|\hbTheta_g^s-\bTheta^*\bigr\|_{\max} 
\le 2K_{\bH^*} \biggl(1+\frac{{ \kappa}}4 \biggr)\lambda. 
\]
\item[(b)] Graphical model selection consistency: picking a
regularization parameter~$\lambda$ to satisfy that
\[
\lambda < \min \biggl( \frac{\psi_{\min}}{2(1+{{\kappa}}/4)K_{\bH^*}},\frac
{d^{-1}}{6(1+{\kappa}/4)K_{\bSigma^*}K_{\bH^*}
\cdot\max\{1,(1+{4}/{\kappa})K_{\bSigma^*}^2K_{\bH^*}\}} \biggr),
\]
then with probability at least $1-p^2\exp(-\frac{\kappa^2}{16}c_0
n\lambda^2)$, $\hbTheta_g^s$ is sign consistent satisfying that
$\sign(\hat\theta^{s.g}_{ij})=\sign(\theta_{ij}^*)$ for any $(i,j)\in
\aset$ and $\hat\theta^{s.g}_{ij}=0$ for any $(i,j)\in\aset^c$.
\end{longlist}
\end{theorem}

In Theorem~\ref{nonasymrglasso}, the condition $\|\bH^*_{\aset^c \aset
}(\bH^*_{\aset\aset})^{-1}\|_{\ell_{\infty}}<1-\kappa$ is also referred
as the \textit{irrepresentable condition}
for studying the theoretical properties of the graphical lasso [\citet
{ravikumar2008}]. We can obtain a straightforward understanding of
Theorem~\ref{nonasymrglasso} by considering its
asymptotic consequences.

\begin{corollary} \label{asymrglasso} Assume that there is a constant
$\kappa\in(0,1)$ such that $\|\bH^*_{\aset^c \aset}(\bH^*_{\aset\aset
})^{-1}\|_{\ell_{\infty}}<1-\kappa$. Suppose that $K_{\bSigma^*}$ and
$K_{\bH^*}$ are both fixed constants.
\begin{longlist}[(a)]
\item[(a)] Rates of convergence: assume $n\gg d^2\log p$, and pick a
regularization parameter $\lambda$ such that $
d^{-1}\gg\lambda=O(({\log p}/n)^{1/2})$. Then we have
\[
\bigl\|\hbTheta_g^s-\bTheta^*\bigr\|_{\max}=O_P
\biggl(\sqrt{\frac{\log p}n} \biggr).
\]
Furthermore, the convergence rates in both Frobenius and matrix $\ell_1$-norms can also be obtained as follows:
\begin{eqnarray*}
\bigl\|\hbTheta_g^s-\bTheta^*
\bigr\|_{F}&=&O_P \biggl(\sqrt{\frac{(s+p)\log p}n} \biggr),
\\
\bigl\|\hbTheta_g^s-\bTheta^*
\bigr\|_{\ell_{1}}&=&O_P \biggl(\sqrt{\frac{\min\{s+p,d^2\}
\log p}n} \biggr).
\end{eqnarray*}
\item[(b)] Graphical model selection consistency: assume $\psi_{\min}$
is also fixed and $n\gg d^2\log p$. Pick a $\lambda$ such that
$d^{-1}\gg\lambda=O(({\log p}/n)^{1/2}).$ Then we have $\sign(\hat\theta^{s.g}_{ij})=\sign(\theta_{ij}^*)$, $\forall(i,j)\in\aset$ and $\sign
(\hat\theta^{s.g}_{ij})=0$, $\forall(i,j)\in\aset^c$.
\end{longlist}
\end{corollary}

Under the same conditions of Theorem~\ref{nonasymrglasso} and
Corollary~\ref{asymrglasso}, by the results in \citet{ravikumar2008},
we know that the conclusions of Theorem~\ref{nonasymrglasso} and
Corollary~\ref{asymrglasso} hold for the oracle graphical lasso. In
other words, the rank-based graphical lasso estimator is comparable to
its oracle counterpart in terms of rates of convergence.

\subsection{\texorpdfstring{Rank-based neighborhood Dantzig selector.}{Rank-based neighborhood Dantzig selector}}\label{sec3.2}

We define $b=\min_k\theta^*_{kk}$, $B=\lambda_{\mathrm{max}}(\bTheta^*)$ and
$M=\|\bTheta^*\|_{\ell_1}$. For each variable $X_k$, define the
corresponding active set $\subaset=\{j\neq k\dvtx \theta^*_{kj}\neq0\}$
with the\vadjust{\goodbreak} maximal cardinality $d=\max_k|\subaset|$. Then we can organize
$\btheta^*_{(k)}$ and $\bTheta_{(k)}^*$ with respect to $\subaset$ as
$\btheta^*_{(k)}=(\btheta^*_{\subaset},\btheta^*_{\subaset^c})$ and
\[
\bTheta_{(k)}^*= \pmatrix{
\bTheta^*_{{\subaset\subaset}} 
&\bTheta^*_{{\subaset\subaset^c}}
\vspace*{2pt}\cr
%
\bTheta^*_{{\subaset^c\subaset}} 
&
\bTheta^*_{{\subaset^c\subaset^c}}}.
\]
Likewise we can partition $\bsigma^*_{(k)}$ and $\bSigma_{(k)}^*$ with
respect to $\subaset$.

\begin{theorem}\label{nonasymds1}
Pick the $\lambda$ such that $d\lambda=o(1)$ and $bn\lambda\ge12 \pi
M$. With probability at least $1-p^2\exp(-c_0\frac{b^2}{M^2}n\lambda^2)$,
there exists $C_{b,B,M}>0$ depending on $b$, $B$ and $M$ only
such that
\[
\bigl\|\bbTheta_{nd}^s-\bTheta^*\bigr\|_{\ell_1}\le 2\bigl\|
\hbTheta_{nd}^s-\bTheta^*\bigr\|_{\ell_1}\le
C_{b,B,M} d\lambda.
\]
\end{theorem}

\begin{corollary}\label{asymds}
Suppose that $b$, $B$ and $M$ are all fixed. Let $n\gg d^2\log p$, and
pick $\lambda$ such that $d^{-1}\gg\lambda=O(({\log p}/n)^{1/2})$. Then
we have
\[
\bigl\|\bbTheta_{nd}^s-\bTheta^*\bigr\|_{\ell_1}=O_P
\biggl( d\sqrt{\frac{\log
p}n} \biggr)\quad\mbox{and}\quad \bigl\|\hbTheta_{nd}^s-
\bTheta^*\bigr\|_{\ell_1}= O_P \biggl( d\sqrt{\frac{\log p}n}
\biggr).
\]
\end{corollary}

\citet{yuan2010} established the rates of convergence of the
neighborhood { Dantzig} selector under the $\ell_1$ norm, which can be
directly applied to the oracle neighborhood Dantzig selector under the
nonparanormal model. Comparing Theorem~\ref{nonasymds1} and
Corollary~\ref{asymds} to the results in \citet{yuan2010}, we see that
the rank-based neighborhood Dantzig selector and the oracle
neighborhood Dantzig selector achieve the same rates of convergence.

Dantzig selector and the lasso are closely related [\citet
{dasso1,dasso2009}]. Similarly to the lasso, the {Dantzig} selector
tends to over-select. \citet{zou2006} proposed the adaptive weighting
idea to develop the adaptive lasso which improves the selection
performance of the lasso and corrects its bias too. The very same idea
can be used to improve the selection performance of Dantzig selector
which leads to the adaptive Dantzig selector [\citet{dicker2009}]. We
can extend the rank-based Dantzig selector to the rank-based adaptive
Dantzig selector. Given adaptive weights $\bw_k$, consider
\begin{equation}
\label{adadantzig}
\quad\hbeta_k^{s.nad}=\mathop{\arg\min}_{\bbeta\in\mathbb{R}^{p-1}}
\| \bw_k\circ\bbeta \|_{\ell_1}\qquad \mbox{subject to } \bigl|
\hbR_{(k)}^s\bbeta-\hbr_{(k)}^s\bigr|\le
\lambda\bw_k,
\end{equation}
where $\circ$ denotes the Hadamard product, and $\mathbf
{a}_{d\times1}\le\mathbf{b}_{d\times1}$ denotes the set of
entrywise inequalities $a_i\le b_i$ for ease of notation. In both our
theoretical analysis and numerical implementation, we utilize the
optimal solution $\hbeta_k^{s.nd}$ of the rank-based Dantzig selector
to construct the adaptive weights $\bw_k$ by
\begin{equation}
\label{adsweight} \bw^d_k= \biggl(\bigl|
\hbeta_k^{s.nd}\bigr|+ \frac{1}n \biggr)^{-1}.
\end{equation}

Define
$
\bbeta^{*}_{\subaset}=(\bTheta^*_{\subaset\subaset})^{-1}\btheta^*_{\subaset},
$
and let
$
\bbeta_k^{*}=(\bbeta^{*}_{\subaset}, \bzero).
$
Thus the support of $\bbeta_k^{*}$ exactly coincides with that of
$\btheta^*_{(k)}$, and then it is further equivalent to the active set
$\subaset$. Define $\psi_k=\|\bbeta^*_{\subaset}\|_{\min}$, $G_k=\|
(\bSigma^*_{{\subaset\subaset}})^{-1}\|_{{ \ell_\infty}}$ and $H_k=\|
\bSigma_{\subaset^c\subaset}^*(\bSigma^*_{{\subaset\subaset}})^{-1}\|_{\ell_\infty}$
for $k=1,2,\ldots,p$. Let $C_0=4B^2(2+\frac{b}M)$.

\begin{theorem}\label{nonasymads}
For each $k$, we pick $\lambda=\lambda_{d}$ as in \eqref{sec22eq2}
satisfying that $\lambda_{d}\ge\frac{12\pi M}{bn}$ and $o(1)=d\lambda_{d}\le\min\{\frac{\psi_{k}}{2C_0}, \frac{1}{4C_0d(\psi_k+2G_k)}-\frac
{1}{C_0n}\}$, and pick $\lambda=\lambda_{ad}$ as in \eqref{adadantzig}
such that
$
\frac{\psi_{k}^2}{8G_k}\ge\lambda_{ad}\ge\max\{\frac{12 \pi}n,
(C_0 d\lambda_{d}+\frac{1}n)\frac{H_k\psi_{k}}{G_k}\},
$
and
$
o(1)=d\lambda_{ad}\le\min\{\lambda_{\min}(\bSigma^*_{\subaset\subaset}),
 \frac{1}{2G_k}, \frac{\psi_{k}}{8G_k(\psi_k+G_k)}\}.
$ In addition, we also choose $\bw_k=\bw_k^d$ as in \eqref{adsweight}
for each $k$. Then with a probability at least $1-p^2\exp(-c_0n\cdot\min
\{\lambda_{ad}^2, \frac{b^2}{M^2}\lambda_{d}^2\})$, for each $k$, the
rank-based adaptive Dantzig selector finds the unique solution $\hbeta^{s.nad}_k=(\hbeta^{s.nad}_{\subaset},\hbeta^{s.nad}_{\subaset^c})$
with $\sign(\hbeta^{s.nad}_{\subaset})=\sign(\bbeta^*_{\subaset})$ and
$\hbeta^{s.nad}_{\subaset^c}=\bzero$, and thus the rank-based
neighborhood adaptive Dantzig selector is consistent for the graphical
model selection.
\end{theorem}


\begin{corollary}\label{asymads}
Suppose $b$, $B$, $M$, $\psi_k$, $G_k$ and $H_k$ ($1\le k\le p$) are
all constants. Assume that $n\gg d^4\log p$ and $\lambda_{\min}(\bSigma^*_{\subaset\subaset})\gg d^2({\log p}/n)^{1/2}$. Pick the tuning
parameters $\lambda_{d}$ and $\lambda_{ad}$ such that $\frac{1}d\gg
\lambda_{d}= O( ({\log p}/n)^{1/2})$ and $\min\{\frac{1}d\cdot\lambda_{\min}(\bSigma^*_{\subaset\subaset}), \frac{1}d\}\gg\lambda_{ad}\gg
d\lambda_{d}$. Then with probability tending to $1$, for each $k$, the
rank-based adaptive Dantzig selector with $\bw_k=\bw_k^d$ as in \eqref
{adsweight} finds the unique optimal solution $\hbeta^{s.nad}_k=(\hbeta^{s.nad}_{\subaset},\hbeta^{s.nad}_{\subaset^c})$ with $\sign(\hbeta^{s.nad}_{\subaset})=\sign(\bbeta^*_{\subaset})$ and $\hbeta^{s.nad}_{\subaset^c}=\bzero$, and thus the rank-based neighborhood
adaptive Dantzig selector is consistent for the graphical model selection.
\end{corollary}

The sign-consistency of the adaptive Dantzig selector is similar to
that of the adaptive lasso [\citet{Geer2011}].
Based on Theorem~\ref{nonasymds1} we construct the adaptive weights in \eqref
{adsweight} which is critical for the success of the rank-based
adaptive Dantzig selector in the high-dimensional setting. It is
important to mention that the rank-based adaptive Dantzig selector does
not require the strong irrepresentable condition for the rank-based
graphical lasso to have the sparsity recovery property. Our treatment
of the adaptive Dantzig selector is fundamentally different from \citet
{dicker2009}. \citet{dicker2009} focused on the canonical linear
regression model and constructed the adaptive weights as the inverse of
the absolute values of ordinary least square estimator. Their
theoretical results only hold in the classical fixed $p$ setting. In
our problem $p$ can be much bigger than $n$. The choice of adaptive
weights in \eqref{adsweight} plays a critical role in establishing the
graphical model selection consistency for the adaptive Dantzig selector
under the high-dimensional setting where $p$ is at a nearly exponential
rate to~$n$. Our technical analysis uses some key ideas such as the
strong duality and the complementary slackness from the linear
optimization theory [\citet{bertsimas1997,boyd2004}].

\subsection{\texorpdfstring{Rank-based CLIME.}{Rank-based CLIME}}\label{sec3.3}
Compared to the graphical lasso, the CLIME can enjoy nice theoretical
properties without assuming the irrepresentable condition [\citet{clime11}].
This continues to hold when comparing the rank-based graphical lasso
and the rank-based CLIME.

\begin{theorem}\label{nonasymclime}
{ Recall that $M=\|\bTheta^*\|_{\ell_1}$.} Pick a regularizing
parameter $\lambda$ such that $n\lambda\ge12 \pi M$. With a
probability at least $1-p^2\exp(-\frac{c_0}{M^2}n\lambda^2)$,
\[
\bigl\|\hbTheta_c^s-\bTheta^*\bigl\|_{\max}\le2M\lambda.
\]
Moreover, assume that $n\gg d^2\log p$, and suppose $M$ is a fixed
constant. Pick a regularization parameter $\lambda$ satisfying $\lambda
=O(({\log p}/n)^{1/2})$. Then we have
\[
\bigl\|\hbTheta_c^s-\bTheta^*\bigr\|_{\max}=O_P
\biggl(\sqrt{\frac{\log p}n} \biggr).
\]
\end{theorem}

Theorem~\ref{nonasymclime} is parallel to Theorem 6 in \citet{clime11}
which can be used to establish the rate of convergence of the oracle CLIME.

To improve graphical model selection performance, \citet{clime11}
suggested an additional thresholding step by applying the element-wise
hard-thresholding rule to $\hbTheta_{c}^s$,
\begin{equation}
\operatorname{HT} \bigl(\hbTheta_{c}^s \bigr)= \bigl(\hat
\theta^{s.c}_{ij}\cdot I_{\{|\hat\theta
^{s.c}_{ij}|\ge\tau_n\}} \bigr)_{1\le i,j\le p},
\end{equation}
where $\tau_n\ge2M\lambda$ is the threshold, and $\lambda$ is given in
Theorem~\ref{nonasymclime}. Here we show that consistent graphical
model selection can be achieved by an adaptive version of the
rank-based CLIME. Given an adaptive weight matrix $\bW$ we define the
rank-based adaptive CLIME as follows:
\begin{equation}
\label{adaptiveclime} {\hbTheta_{ac}^s}=\mathop{\arg\min}_{\bTheta}\|\bW\circ\bTheta\|_1 \qquad\mbox{subject to } \bigl|
\hbR^s\bTheta-\bI\bigr|\le\lambda\bW,
\end{equation}
where $\mathbf{A}_{p\times p}\le\mathbf{B}_{p\times p}$ is a
simplified expression for the set of inequalities $a_{ij}\le b_{ij}$
(for all $1\le i,j\le p$).
Write $\bW=(\bw_1,\ldots,\bw_p)$. By Lemma 1 in \citet{clime11} the
above linear programming problem in (\ref{adaptiveclime}) is exactly
equivalent to $p$ vector minimization subproblems,
\[
{\hbtheta^{s.ac}_k}=\mathop{\arg\min}_{\btheta\in\mathbb{R}^p}\|
\bw_k\circ\btheta \|_{\ell_1}\qquad \mbox{subject to } \bigl|
\hbR^s\btheta-\bolde_k\bigr|\le \lambda\bw_k.
\]
In both our theory and implementation, we utilize the rank-based
CLIME's optimal solution $\hbTheta_c^s$ to construct an adaptive weight
matrix $\bW$ by
\begin{equation}
\label{aclimeweight} \bW^c= \biggl(\bigl|\hbTheta_c^s\bigr|+
\frac{1}n \biggr)^{-1}.
\end{equation}

We now prove the graphical model selection consistency\vspace*{1pt} of the
rank-based adaptive CLIME. Denote $\bTheta^*$ as $(\btheta_1^*,\ldots
,\btheta_p^*)$, and define $\tilde\subaset=\subaset\cup\{k\}$. Then we
can organize $\btheta^*_k$ and $\bSigma^*$ with respect to $\tilde
\subaset$ and $\tilde{\subaset^c}$. For $k=1,2,\ldots,p$, define
$\tilde G_k=\|(\bSigma^*_{{\tilde\subaset\tilde\subaset}})^{-1}\|_{{
\ell_\infty}}$ and $\tilde H_k=\|\bSigma_{\tilde{\subaset^c}\tilde
\subaset}^*(\bSigma^*_{{\tilde\subaset\tilde\subaset}})^{-1}\|_{\ell
_\infty}$.

\begin{theorem}\label{nonasymaclime}
Recall $\psi_{\min}=\min_{(i,j)\in\aset}|\theta^*_{ij}|$. For each
$k$ pick $\lambda=\lambda_{c}$ as in~\eqref{sec22eq3} such that $\min\{
\frac{\psi_{\min}}{4M}, \frac{1}{4M(\psi_{\min}+2\tilde G_k)d}-\frac
{2}{Mn}\}
\ge\lambda_{c}\ge\frac{12\pi M}n$ and $d\lambda_{c}=o(1)$, and we
further pick the regularization parameter $\lambda=\lambda_{ac}$ as in
\eqref{adaptiveclime} satisfying that $\frac{\psi^2_{\min}}{8\tilde G_k}
\ge\lambda_{ac}
\ge\max\{{12 \pi}/n, (2M\lambda_{c}+\frac{1}n)\frac{\tilde H_k\psi_{\min}}{\tilde G_k}\}$ and
$o(1)=d\lambda_{ac}\le\min\{\lambda_{\min
}(\bSigma^*_{\subaset\subaset}), \frac{1}{2\tilde G_k}, \frac{\psi_{\min
}}{4\tilde G_k(\psi_{\min}+\tilde G_k)}\}
$. In addition we choose $\bW=\bW^c$ as in \eqref{aclimeweight}. With
a probability at least $1-p^2\exp(-c_0n\min\{\lambda_{ac}^2, \frac
{1}{M^2}\lambda_{c}^2\})$, the rank-based adaptive CLIME's optimal
solution $\hbTheta_{ac}^s$ is sign consistent, that is, $\sign(\hat
\theta^{s.ac}_{ij})=\sign(\theta_{ij}^*)$ for $(i,j)\in\aset$ and $\sign
(\hat\theta^{s.ac}_{ij})=0$ for $(i,j)\in\aset^c$.
\end{theorem}

\begin{corollary}\label{asymaclime}
Suppose $\psi_{\min}$, $M$, $\tilde G_k$ and $\tilde H_k$ ($1\le k\le
p$) are all constants. Assume that $n\gg d^2\log p$ and $\lambda_{\min
}(\bSigma^*_{\tilde\subaset\tilde\subaset})\gg d{({\log
p}/n)^{1/2}}$. Pick the regularization parameters $\lambda_{c}$ and
$\lambda_{ac}$ such that $\frac{1}d\ge\lambda_{c}= O( ({\log
p}/n)^{1/2}), $ and $\min\{\lambda_{\min}(\bSigma^*_{\tilde\subaset
\tilde\subaset})/d, \frac{1}d\}\gg\lambda_{ac}\gg\lambda_{c}$. Let $\bW
=\bW^c$ as in \eqref{aclimeweight}. Then with probability tending to
$1$, the rank-based adaptive CLIME's optimal solution $\hbTheta_{ac}^s$
is sign consistent for the graphical model selection, that is, $\sign
(\hat\theta^{s.ac}_{ij})=\sign(\theta_{ij}^*)$ for $(i,j)\in\aset$ and
$\sign(\hat\theta^{s.ac}_{ij})=0$ for $(i,j)\in\aset^c$.
\end{corollary}

The nice theoretical property of the rank-based CLIME allows us to
construct the adaptive weights in \eqref{aclimeweight}, which is
critical for establishing the graphical model selection consistency for
the rank-based adaptive CLIME estimator in the high-dimensional setting
without the strong ir-representable condition.

\section{\texorpdfstring{Numerical properties.}{Numerical properties}}\label{sec4}
In this section we present both simulation studies and real examples to
demonstrate the finite sample performance of the proposed rank-based estimators.

\subsection{\texorpdfstring{Monte Carlo simulations.}{Monte Carlo simulations}}\label{sec4.1}

In the simulation study, we consider both Gaussian data and
nonparanormal data. In models 1--4 we draw $n$ independent samples from
$N_p(0,\bSigma)$ with four
different $\bTheta$:
\begin{longlist}[Model 1:]
\item[\textit{Model} 1:] $\theta_{ii}=1$ and $\theta_{i,i+1}=0.5$;
\item[\textit{Model} 2:] $\theta_{ii}=1$, $\theta_{i,i+1}=0.4$ and $\theta_{i,i+2}=\theta_{i,i+3}=0.2$;
\item[\textit{Model} 3:] Randomly choose $16$ nodes to be the hub nodes in
$\bTheta$, and each of them connects with
$5$ distinct\vadjust{\goodbreak} nodes with $\bTheta_{ij}=0.2$. Elements, not associated
with hub nodes, are set as $0$
in $\bTheta$. The diagonal element $\sigma$ is chosen similarly as
that in the previous model.
\item[\textit{Model} 4:] $\bTheta=\bTheta_0+\sigma I$, where $\bTheta_0$ is a
zero-diagonal symmetric matrix. Each off-diagonal element ${\bTheta_0}_{ij}$ independently follows
a point mass $0.99\delta_0+0.01\delta_{0.2}$, and the diagonal element $\sigma$ is set to be the absolute
value of the minimal negative eigenvalue of $\bTheta_0$ to ensure the
semi-positive-definiteness of $\bTheta$.
\end{longlist}

In models 1b--4b we first generate $n$ independent data from
$N_p(0,\bSigma)$ and then transfer the normal data using transformation
functions
\[
\mathbf{g}= \bigl[f^{-1}_1, f^{-1}_2,
f^{-1}_3, f^{-1}_4,
f^{-1}_5, f^{-1}_1,
f^{-1}_2, f^{-1}_3,
f^{-1}_4, f^{-1}_5, \ldots \bigr],
\]
where $f_1(x)=x$, $f_2(x)=\log(x)$, $f_3(x)=x^{\frac{1}3},$ $f_4(x)=\log
(\frac{x}{1-x})$ and $f_5(x)=f_2(x)I_{\{x<-1\}}+f_1(x)I_{\{-1\le x\le
1\}}+(f_4(x-1)+1)I_{\{x>1\}}$.
In all cases we let $n=300$ and $p=100$.

%

\begin{table}
\caption{List of all estimators in the numercial study}\label{simtab1}
\begin{tabular*}{\textwidth}{@{\extracolsep{\fill}}lp{235pt}@{}}
\hline
\multicolumn{1}{@{}l}{\textbf{Notation}} & \textbf{Details} \\
\hline
GLASSO & Penalized likelihood estimation via graphical lasso\\
MB & Neighborhood lasso [\citet{meinshausen2006}]\\
MB.au (or MB.ai) & MB${}+{}$aggregation by union (or by intersection)\\
NDS & Neighborhood selection via Dantzig selector \\
NDS.au (or NDS.ai) & NDS${}+{}$aggregation by union (or by intersection)\\
CLIME & Constrained $\ell_1$ minimization estimator \\
LLW & The ``plug-in'' extension of GLASSO [\citet{lafferty2009}]\\
R-GLASSO& Proposed rank-based extension of GLASSO \\
R-NDS & Proposed rank-based extension of NDS\\
R-NDS.au (or R-NDS.ai) & R-NDS${}+{}$aggregation by union (or by
intersection)\\
R-NADS & Proposed rank-based {adaptive extension of R-NDS}\\
R-NADS.au (or R-NADS.ai) & R-NADS${}+{}$aggregation by union (or by
intersection)\\
R-CLIME & Proposed rank-based extension of CLIME \\
R-ACLIME& Proposed rank-based adaptive extension of CLIME \\
\hline
\end{tabular*}
\end{table}
%

\begin{table}
\caption{Estimation performance in the Gaussian model. Estimation
accuracy is measured by the matrix $\ell_2$-norm with standard errors
in the bracket}\label{simtab2}
\begin{tabular*}{\textwidth}{@{\extracolsep{\fill}}lcccc@{}}
\hline
\textbf{Method} & \textbf{Model 1} & \textbf{Model 2} & \textbf{Model 3} & \multicolumn{1}{c@{}}{\textbf{Model 4}} \\
\hline
{GLASSO}
&0.74 &1.23 &0.67 &0.63 \\
&(0.01) &(0.02) &(0.01) &(0.01) \\
{LLW}
&0.84 &1.28 &0.68 &0.67 \\
&(0.01) &(0.02) &(0.01) &(0.01) \\
{R-GLASSO}
&0.81 &1.30 &0.64 &0.70 \\
&(0.01) &(0.02) &(0.01) &(0.01) \\[3pt]
{NDS}
&0.78 &1.25 &0.61 &0.57 \\
&(0.01) &(0.02) &(0.01) &(0.01) \\
{R-NDS}
&0.81 &1.28 &0.63 &0.62 \\
&(0.01) &(0.02) &(0.01) &(0.01) \\[3pt]
{CLIME}
&0.71 &1.19 &0.54 &0.59 \\
&(0.01) &(0.02) &(0.01) &(0.01) \\
{R-CLIME}
&0.79 &1.27 &0.58 &0.61 \\
&(0.01) &(0.02) &(0.01) &(0.01) \\
\hline
\end{tabular*}
\end{table}

\begin{table}[b]
\caption{Estimation performance in the nonparanormal model. Estimation
accuracy is measured by the matrix $\ell_2$-norm with standard errors
in the bracket}\label{simtab2b}
\begin{tabular*}{\textwidth}{@{\extracolsep{\fill}}lcccc@{}}
\hline
\textbf{Method} & \textbf{Model 1b} & \textbf{Model 2b} & \textbf{Model 3b} & \multicolumn{1}{c@{}}{\textbf{Model 4b}}\\
\hline
{GLASSO}
&1.77 &2.68 &1.31 &1.28 \\
&(0.01) &(0.06) &(0.02) &(0.01) \\
\multirow{2}{*}{LLW}
&0.84 &1.28 &0.68 &0.67 \\
&(0.01) &(0.01) &(0.01) &(0.01) \\
{R-GLASSO}
&0.81 &1.30 &0.64 &0.70 \\
&(0.01) &(0.02) &(0.01) &(0.01) \\[3pt]
{NDS}
&1.41 &2.42 &1.16 &1.13 \\
&(0.01) &(0.03) &(0.02) &(0.02) \\
{R-NDS}
&0.81 &1.28 &0.63 &0.62 \\
&(0.01) &(0.02) &(0.01) &(0.01) \\[3pt]
{CLIME}
&1.22 &2.51 &1.24 &1.03 \\
&(0.02) &(0.03) &(0.02) &(0.01) \\
{R-CLIME}
&0.79 &1.27 &0.58 &0.61 \\
&(0.01) &(0.02) &(0.01) &(0.01) \\
\hline
\end{tabular*}
\end{table}

Table~\ref{simtab1} summarizes all the estimators investigated in our study.
For each estimator, the tuning parameter is chosen by cross-validation.
Estimation accuracy is measured by the average {matrix $\ell_2$-norm}
over 100 independent replications, and selection accuracy is evaluated
by the average false positive/negative.

\begin{table}
\caption{Selection performance in the Gaussian model. Selection
accuracy is measured by counts of false negative (\#FN) or false
positive (\#FP) with standard errors in the bracket}\label{simtab3}
\begin{tabular*}{\textwidth}{@{\extracolsep{\fill}}lcccccccc@{}}
\hline
& \multicolumn{2}{c}{\textbf{Model 1}} & \multicolumn{2}{c}{\textbf{Model 2}} &
\multicolumn{2}{c}{\textbf{Model 3}}
& \multicolumn{2}{c}{\textbf{Model 4}} \\[-6pt]
& \multicolumn{2}{c}{\hrulefill} & \multicolumn{2}{c}{\hrulefill} &
\multicolumn{2}{c}{\hrulefill}
& \multicolumn{2}{c@{}}{\hrulefill} \\
& \textbf{\#FN} & \textbf{\#FP} & \textbf{\#FN} & \textbf{\#FP} & \textbf{\#FN} & \textbf{\#FP} & \textbf{\#FN} & \textbf{\#FP} \\
\hline
{GLASSO}
&0.00 &521.21 &263.16 &45.21 &0.00 &114.48 &0.03 &35.33 \\
&(0.00) &(1.91) &(0.58) &(1.26) &(0.00) &(1.94) &(0.02) &(1.29) \\
{LLW}
&0.00 &518.84 &264.18 &43.45 &0.00 &116.02 &0.04 &35.08 \\
&(0.00) &(1.91) &(0.56) &(1.34) &(0.00) &(2.01) &(0.02) &(1.19) \\
{R-GLASSO}
&0.00 &505.77 &264.86 &48.01 &0.00 &114.89 &0.03 &37.13 \\
&(0.00) &(1.67) &(0.57) &(1.57) &(0.00) &(2.17) &(0.02) &(1.07) \\[3pt]
{MB.au}
&0.00 &154.81 &232.99 &89.61 &0.00 &44.03 &0.02 &41.22 \\
&(0.00) &(1.29) &(0.74) &(1.37) &(0.00) &(0.81) &(0.01) &(0.77) \\
{R-NDS.au}
&0.00 &163.78 &230.77 &118.46 &0.00 &69.16 &0.03 &49.31 \\
&(0.00) &(1.27) &(0.79) &(2.12) &(0.00) &(0.92) &(0.02) &(0.88) \\
{R-NADS.au}
&0.00 &80.90 &218.69 &83.62 &0.00 &60.75 &0.03 &48.59 \\
&(0.00) &(2.52) &(1.02) &(2.90) &(0.00) &(1.04) &(0.02) &(0.92) \\[3pt]
{MB.ai}
&0.00 &30.62 &260.76 &21.79 &0.00 &9.42 &0.04 &9.58 \\
&(0.00) &(0.53) &(0.55) &(0.60) &(0.00) &(0.31) &(0.02) &(0.34) \\
{R-NDS.ai}
&0.00 &38.62 &259.66 &29.34 &0.00 &11.52 &0.07 &11.87 \\
&(0.00) &(0.52) &(0.61) &(0.68) &(0.00) &(0.40) &(0.04) &(0.40) \\
{R-NADS.ai}
&0.06 &14.92 &256.16 &24.62 &0.00 &10.54 &0.08 &10.98 \\
&(0.02) &(0.11) &(0.68) &(0.79) &(0.00) &(0.36) &(0.04) &(0.38) \\[3pt]
{CLIME}
&0.00 &143.88 &263.77 &34.71 &0.00 &32.53 &0.02 &32.59 \\
&(0.00) &(0.10) &(0.57) &(1.42) &(0.00) &(0.78) &(0.01) &(1.17) \\
{R-CLIME}
&0.00 &148.24 &265.81 &38.23 &0.00 &37.44 &0.04 &36.56 \\
&(0.01) &(3.11) &(1.22) &(2.55) &(0.05) &(2.45) &(0.33) &(1.18) \\
{R-ACLIME}
&0.00 &82.53 &264.74 &34.52 &0.00 &29.83 &0.07 &31.09 \\
&(0.00) &(0.13) &(0.63) &(2.60) &(0.00) &(0.61) &(0.03) &(1.02) \\
\hline
\end{tabular*}
\end{table}

\begin{table}
\caption{Selection performance in the nonparanormal model. Selection
accuracy is~measured~by~counts of false negative (\#FN) or false
positive (\#FP) with~standard~errors in the bracket}\label{simtab4}

\begin{tabular*}{\textwidth}{@{\extracolsep{\fill}}lcccccccc@{}}
\hline
& \multicolumn{2}{c}{\textbf{Model 1b}} & \multicolumn{2}{c}{\textbf{Model 2b}} &
\multicolumn{2}{c}{\textbf{Model 3b}}
& \multicolumn{2}{c@{}}{\textbf{Model 4b}} \\[-4pt]
& \multicolumn{2}{c}{\hrulefill} & \multicolumn{2}{c}{\hrulefill} &
\multicolumn{2}{c}{\hrulefill}
& \multicolumn{2}{c@{}}{\hrulefill} \\
& \textbf{\#FN} & \textbf{\#FP} & \textbf{\#FN} & \textbf{\#FP} & \textbf{\#FN} & \textbf{\#FP} & \textbf{\#FN} & \textbf{\#FP} \\
\hline
{GLASSO}
&58.81 &470.05 &286.40 &44.70 &9.82 &134.70 &8.06 &44.20 \\
&(0.35) &(5.30) &(0.74) &(1.48) &(0.41) &(2.08) &(0.36) &(1.33) \\
{LLW}
&0.00 &518.84 &264.18 &43.45 &0.00 &116.02 &0.04 &35.08 \\
&(0.00) &(1.91) &(0.56) &(1.34) &(0.00) &(2.01) &(0.02) &(1.19) \\
{R-GLASSO}
&0.00 &505.77 &264.86 &48.01 &0.00 &114.89 &0.03 &37.13 \\
&(0.00) &(1.67) &(0.57) &(1.57) &(0.00) &(2.17) &(0.02) &(1.07) \\[3pt]
{MB.au}
&56.28 &472.86 &283.15 &61.69 &12.99 &99.10 &8.28 &57.65 \\
&(0.26) &(4.11) &(0.64) &(1.04) &(0.46) &(1.31) &(0.36) &(0.90) \\
{R-NDS.au}
&0.00 &163.78 &230.77 &118.46 &0.00 &69.16 &0.03 &49.31 \\
&(0.00) &(1.27) &(0.79) &(2.12) &(0.00) &(0.92) &(0.02) &(0.88) \\
{R-NADS.au}
&0.00 &80.90 &218.69 &83.62 &0.00 &60.75 &0.03 &48.59 \\
&(0.00) &(2.52) &(1.02) &(2.90) &(0.00) &(1.04) &(0.02) &(0.92) \\[3pt]
{MB.ai}
&68.68 &197.44 &304.71 &22.72 &16.88 &50.25 &11.67 &23.88 \\
&(0.16) &(1.12) &(0.61) &(0.56) &(0.52) &(0.92) &(0.42) &(0.50) \\
{R-NDS.ai}
&0.00 &38.62 &259.66 &29.34 &0.00 &11.52 &0.08 &11.87 \\
&(0.00) &(0.52) &(0.61) &(0.68) &(0.00) &(0.40) &(0.04) &(0.40) \\
{R-NADS.ai}
&0.06 &14.92 &256.16 &24.62 &0.00 &10.54 &0.08 &10.98 \\
&(0.02) &(0.11) &(0.68) &(0.79) &(0.00) &(0.36) &(0.04) &(0.38) \\[3pt]
{CLIME}
&47.14 &385.95 &286.16 &45.25 &10.02 &123.31 &7.87 &46.38 \\
&(0.39) &(1.99) &(0.74) &(1.45) &(0.41) &(2.11) &(0.36) &(1.34) \\
{R-CLIME}
&0.00 &148.24 &265.81 &38.23 &0.00 &37.44 &0.04 &36.56 \\
&(0.01) &(3.11) &(1.22) &(2.55) &(0.05) &(2.45) &(0.33) &(1.18) \\
{R-ACLIME}
&0.00 &82.53 &264.74 &34.52 &0.00 &29.83 &0.07 &31.09 \\
&(0.00) &(0.13) &(0.63) &(2.60) &(0.00) &(0.61) &(0.03) &(1.02) \\
\hline
\end{tabular*}
\end{table}

The simulation results are summarized in Tables~\ref{simtab2}--\ref
{simtab4}. First of all, we can see that the graphical lasso,
neighborhood selection and CLIME do not have satisfactory performance
under models 1b--4b due to the lack of ability to handle nonnormality.
Second, the three rank-based estimators perform similarly to their
oracle counterparts. Note that in models 1b--4b the oracle graphical
lasso, the oracle neighborhood Dantzig and the oracle CLIME are
actually the graphical lasso, the neighborhood Dantzig and the CLIME in
models 1--4. In terms of precision matrix estimation the rank-based
CLIME seems to be the best, while the rank-based neighborhood adaptive
Dantzig selector has the best graphical model selection performance. We
have also obtained the simulation results under the matrix $\ell_1$-norm.
The conclusions stay the same. For space consideration we
leave these $\ell_1$-norm results to the technical report version of
this paper.

\subsection{\texorpdfstring{Applications to gene expression genomics.}{Applications to gene expression genomics}}\label{sec4.2}
We illustrate our proposed rank-based estimators on a real data set
to recover the isoprenoid genetic regulatory network in Arabidposis
thaliana [\citet{wille2004}]. This dataset contains
the gene expression measurements of 39 genes (excluding protein GGPPS7
in the MEP pathway) assayed on $n = 118$ Affymetrix GeneChip
microarrays.

We used seven estimators (GLASSO, MB, CLIME, LLW, R-GLASSO, R-NADS and
R-ACLIME) to reconstruct the regulatory network. The first three
estimators are performed after taking the log-transformation of the
original data, and the other four estimators are directly applied to
the original data. To be more conservative, we only considered the
integration by union for the neighborhood selection procedures. We
generated $100$ independent Bootstrap samples and computed the
frequency of each edge being selected by each estimator. The final
model by each method only includes
edges selected by at least $80$ times over $100$ Bootstrap samples. We
report the number of selected edges by each estimator in Table~\ref{rankreal}.
The rank-based graphical lasso performs similarly to the LLW method.
The rank-based adaptive CLIME produces the sparsest graphs. {We also
compared pairwise intersections of the selected edges among different
estimators. More than $70\%$ of the selected edges by GLASSO, MB or
CLIME turn out to be validated by both LLW and R-GLASSO, and more than
$40\%$ of the selected edges by GLASSO, MB or CLIME are justified by
R-NADS and R-ACLIME.} The selected models support the biological
arguments that the interactions between the pathways do exist although
they operate
independently under normal conditions
[\citet{laule2003,rodriguez2004}].

\begin{table}
\caption{The isoprenoid genetic regulatory network: counts of stable
edges}\label{rankreal}
\begin{tabular*}{\textwidth}{@{\extracolsep{\fill}}lcccc@{}}
\hline
& \textbf{GLASSO} & \textbf{Neighborhood LASSO} & \textbf{CLIME} &\\
\hline
$\#$ of stable edges & 100 & 101 & 67 &\\
\hline
& \textbf{LLW} & \textbf{R-GLASSO} & \textbf{R-NADS} & \textbf{R-ACLIME} \\
\hline
$\#$ of stable edges & 87 & 88 & 50 & 52 \\
\hline
\end{tabular*}
\end{table}

\section{\texorpdfstring{Discussion.}{Discussion}}\label{sec5}

Using ranks of the raw data for statistical inference is a powerful and
elegant idea in the nonparametric statistics literature; see \citet
{lehmann1998} for detailed treatment and discussion.
Some classical rank-based statistical methods include Friedman's test
in analysis of variance and Wilcoxon signed-rank test. This work is
devoted to the rank-based estimation of $\bSigma^{-1}$ of the
nonparanormal model under a strong sparsity assumption that $\bSigma^{-1}$ has only a few nonzero entries, and our results show that
rank-based estimation is still powerful and elegant in the new setting
of high-dimensional nonparametric graphical modeling. In a separate
paper, \citet{xue2011} also studied the problem of optimal estimation of
$\bSigma$ of the nonparanormal model under a weak sparsity assumption
that $\bSigma$ belongs to some weak $\ell_q$ ball and showed that a
rank-based thresholding estimator is adaptive minimax optimal under the
matrix $\ell_1$ norm and $\ell_2$ norm.

\begin{appendix}
\section*{\texorpdfstring{Appendix: Technical proofs}{Appendix: Technical
proofs}}\label{app}\vspace*{-3pt}
\begin{pf*}{Proof of Theorem~\ref{nonasymrglasso}}
Using Lemma 3 in \citet{ravikumar2008}, $\hbTheta_g^s\succ0$ is
uniquely characterized by the sub-differential optimality condition that
$\hbR^s-(\hbTheta_g^s)^{-1}+\lambda\hbZ= \bzero,$
where
$
\hbZ
$ is the sub-differential with respect to $\hbTheta_g^s$. Define the
``oracle'' estimator $\tbTheta_g^s$ by
\[
\tbTheta^s_g=\mathop{\arg\min}_{\bTheta\succ0, \bTheta_{\aset^c}=\bzero}-\log \det(
\bTheta)+\tr \bigl(\hbSigma^o\bTheta \bigr)+\lambda\sum
_{i\neq j} |\theta_{ij}|. 
\]
Then we can construct $\tbZ$ to satisfy that $\hbR^s-(\tbTheta^s_g)^{-1}+\lambda\tbZ= \bzero.$
As in \citet{ravikumar2008} the rest of the proof depends on a
exponential-type concentration bound concerning the accuracy of the
sample estimator of the correlation matrix under the entry-wise $\ell_\infty$ bound.
Our Lemma~\ref{concentration} fulfills that role. With
Lemma~\ref{concentration}, Theorem~\ref{nonasymrglasso} can be proved by following the
line of the proof in \citet{ravikumar2008}. For the sake of space, we
move the rest of proof to the technical report version of this paper
[\citet{XZ11b}].

We now prove Lemma~\ref{concentration}. First, Spearman's rank correlation $\hat r_{ij}$
can be written in terms of the Hoeffding decomposition [\citet{hoeffding1948}]
\begin{equation}
\label{hoeff} \hat r_{ij} = \frac{n-2}{n+1}u_{ij}+
\frac{3}{n+1}d_{ij},
\end{equation}
where
$
d_{ij}=\frac{1}{n(n-1)}\sum_{k\neq l}
\sign(x_{ki}-x_{li})
\cdot\sign(x_{kj}-x_{lj}),
$
and
\begin{equation}
\label{Ustat1} u_{ij} = \frac{3}{n(n-1)(n-2)} \sum
_{k\neq l,k\neq m,l\neq m} \sign(x_{ki}-x_{li})\cdot
\sign(x_{kj}-x_{mj}).
\end{equation}
Direct calculation yields that $\E(u_{ij})=\frac{6}{\pi}\sin^{-1}(\frac
{\sigma_{ij}}2)$ [\citet{kendall1948}].
%
%
Then we can obtain that $\sigma_{ij}=2\sin(\frac{\pi}{6}E(u_{ij}))$. By
definition $\hat r_{ij}^s=2\sin(\frac{\pi}{6}\hat r_{ij})$. Note that
$2\sin(\frac{\pi}6\cdot)$ is a Lipschitz function with the Lipschitz
constant $\pi/3$. Then we have
\[
\Pr \bigl(\bigl|\hat r^s_{ij}-\sigma_{ij}\bigr|>
\varepsilon \bigr) \le \Pr \biggl(\bigl|\hat r_{ij}-E(u_{ij})\bigr|>
\frac{3\varepsilon}{\pi} \biggr).
\]
Applying (\ref{hoeff}) and (\ref{Ustat1})\vspace*{1pt} yields
$
\hat r_{ij}-E(u_{ij})=u_{ij}-E(u_{ij})+\frac{3}{n+1}d_{ij}-\frac{3}{n+1}u_{ij}.
$
Note $|u_{ij}|\le3$ and $|d_{ij}|\le1$. Hence, $|u_{ij}|\le\frac
{\varepsilon}{4\pi}(n+1)$ and $|d_{ij}|\le\frac{\varepsilon}{4\pi
}(n+1)$ always hold provided that $n>12\pi/\varepsilon$, which are
satisfied by the assumption in Lemma~\ref{concentration}. For such chosen $n$, we have
\[
\Pr \biggl(\bigl|r_{ij}-E(u_{ij})\bigr|>\frac{3\varepsilon}{\pi} \biggr) \le
\Pr \biggl(\bigl|u_{ij}-E(u_{ij})\bigr|>\frac{3\varepsilon}{2\pi} \biggr).
\]
Finally, we observe that $u_{ij}$ is a function of independent samples
$(\bx_1,\ldots,\bx_n)$. Now we make a claim that if we replace the
$t$th sample by some $\hat{\bx}_{t}$, the change in $u_{ij}$ will be
bounded as
\begin{equation}
\sup_{\bx_1,\ldots,\bx_n,\hat{\bx}_{t}} \bigl|u_{ij}(\bx_1,\ldots,
\bx_n) -u_{ij}(\bx_1,\ldots,
\bx_{t-1},\hat{\bx}_{t},\bx_{t+1},\ldots,
\bx_n)\bigr| \le\frac{15}{n}.\label{step03}
\end{equation}
Then we can apply the McDiarmid's inequality [\citet{mcdiarmid1989}] to
conclude the desired concentration bound for some absolute constant $c_0>0$,
\[
\Pr \bigl(\bigl|\hat r^s_{ij}-\sigma_{ij}\bigr|>
\varepsilon \bigr) \le \Pr \biggl(\bigl|u_{ij}-E(u_{ij})\bigr|\ge
\frac{3\varepsilon}{2\pi} \biggr) \le 
2\exp \bigl(-c_0n
\varepsilon^2 \bigr).
\]

Now it remains to verify (\ref{step03}) to complete the proof of Lemma
\ref{concentration}. We provide a brief proof for this claim. Assume
that $\bx_t=(x_{1t},\ldots,x_{{ p}t})'$ is replaced by $\tilde{\bx
}_t=(\tilde{x}_{1t},\ldots,\tilde{x}_{{ p}t})'$, and we want to prove
that the change of $u_{ij}$ is at most $15/n$. Without loss of
generality we may assume that $n_i=\#\{s\dvtx \sign(\tilde
x_{ti}-x_{si})=-\sign(x_{ti}-x_{si}), s\neq t\}$ and also assume that\vadjust{\goodbreak}
$n_j=\#\{s\dvtx \sign(\tilde x_{tj}-x_{sj})=-\sign(x_{tj}-x_{sj}), s\neq t\}$.
Then we have
\begin{eqnarray*}
&& \bigl|u_{ij}(\bx_1,\ldots,\bx_n)-u_{ij}(
\bx_1,\ldots,\bx_{t-1},\tilde{\bx }_{t},
\bx_{t+1},\ldots,\bx_n) \bigr|
\\
&&\qquad\le \biggl|\sum_{k\neq t,k\neq m,m\neq t} \bigl\{\sign(x_{ki}-x_{ti})-
\sign(x_{ki}-\tilde x_{ti}) \bigr\}\cdot\sign
(x_{kj}-x_{mj})
\\
&&\qquad\quad{}+ \sum_{k\neq t,k\neq l,l\neq t} \bigl(\sign(x_{kj}-x_{tj})-
\sign(x_{kj}-\tilde x_{tj}) \bigr)\cdot\sign
(x_{ki}-x_{li})
\\
&&\qquad\quad{}+ \sum_{l\neq t,m\neq t,l\neq m} \bigl(\sign(x_{ti}-x_{li})
\cdot\sign(x_{tj}-x_{mj})-\sign(\tilde x_{ti}-x_{li})
\\
&&\hspace*{198pt}\qquad\quad{}\times\sign(\tilde x_{tj}-x_{mj}) \bigr)\biggr |\\
&&\quad\qquad{}\times
\frac{3}{n(n-1)(n-2)}
\\
&&\qquad\le \frac{3\cdot2\cdot
[n_i(n-2)+n_j(n-2)+n_j(n-1-n_i)+n_i(n-1-n_j)]}{n(n-1)(n-2)}
\\
&&\qquad\le \frac{12}n \biggl(1+\frac{1}4\frac{1}{(n-1)(n-2)} \biggr)
\\
&&\qquad\le \frac{15}{n},
\end{eqnarray*}
where the third inequality holds if and only if $n_i=n_j=n-\frac
{3}2. $
\end{pf*}

\begin{pf*}{Proof of Theorem~\ref{nonasymds1}} For space of
consideration, we only show the sketch of the proof, and the detailed
proof is relegated to the supplementary file [\citet{supp}] and also the technical
report version of this paper [\citet{XZ11b}]. We begin with an
important observation that we only need to prove the risk bound for
$\hbTheta^s_{nd}$ because
\[
\bigl\|\bbTheta^s_{nd}-\bTheta^*\bigr\|_{\ell_1} 
\le \bigl\|\bbTheta^s_{nd}-\hbTheta^s_{nd}
\bigr\|_{\ell_1} + \bigl\|\hbTheta^s_{nd}-\bTheta^*
\bigr\|_{\ell_1} \le 2\bigl\|\hbTheta^s_{nd}-\bTheta^*
\bigr\|_{\ell_1}.
\]
To bound the difference between $\hbTheta^s_{nd}$ and $\bTheta^*$ under
the matrix $\ell_1$-norm, we only need to bound $|\hat\theta^{s.nd}_{kk}-\theta^*_{kk}|$ and
$\|\hbtheta^{s.nd}_{(k)}-\btheta^*_{(k)}\|_{\ell_1}$ for each $k=1,\ldots,p$. To this end, we consider
the probability event $\{\|\hbR^s-\bSigma^*\|_{\max}\le\frac{b}M\lambda
\}$, and under this event, we can show that for $k=1,\ldots,p$,
\begin{equation}
\bigl\|\hbR^s_{(k)}\bbeta_k^*-\hbr^s_{(k)}
\bigr\|_{\ell_{\infty}}\le\lambda\quad \mbox{and}\quad \bigl\|\hbeta_k^{s.nd}-
\bbeta_k^*\bigl\|_{\ell_1}\le C_0 d\lambda,
\label{rankdantzigstep0}
\end{equation}
where $C_0$ is some quantity depending on $b$, $B$ and $M$ only.

Now we can use \eqref{rankdantzigstep0} to further bound $|\hat\theta^{s.nd}_{kk}-\theta^*_{kk}|$
under the same event. To this end, we
first derive an upper bound for $|(\hat\theta^{s.nd}_{kk})^{-1}-(\theta_{kk}^*)^{-1}|$ as
\begin{equation}
\bigl| \bigl(\hat\theta^{s.nd}_{kk} \bigr)^{-1}- \bigl(
\theta_{kk}^* \bigr)^{-1}\bigr| \le \biggl(1+\frac{b}M
\biggr)\cdot\lambda\bigl\|\bbeta_k^*\bigr\|_{\ell_1} + \bigl\|
\hbeta_{k}^{s.nd}-\bbeta_{k}^*\bigr\|_{\ell_1}.
\label{rankdantzigstep1}
\end{equation}
Notice that $|\hat\theta^{s.nd}_{kk}-\theta^*_{kk}|
=
|(\hat\theta^{s.nd}_{kk})^{-1}-(\theta_{kk}^*)^{-1}|\cdot|\hat\theta^{s.nd}_{kk}|\cdot|\theta^*_{kk}|$ and also
$
|\hat\theta^{s.nd}_{kk}|\le|\hat\theta^{s.nd}_{kk}-\theta^*_{kk}|+|\theta^*_{kk}|
$. Then $|\hat\theta^{s.nd}_{kk}-\theta^*_{kk}|$ can be upper bounded by
\begin{eqnarray*}
\bigl|\hat\theta^{s.nd}_{kk}-\theta^*_{kk}\bigr| &\le&
\frac{|(\hat\theta^{s.nd}_{kk})^{-1}-(\theta_{kk}^*)^{-1}|\cdot|\theta^*_{kk}|^2} {
1-|(\hat\theta^{s.nd}_{kk})^{-1}-(\theta_{kk}^*)^{-1}|\cdot|\theta^*_{kk}|}\\[-1pt]
& \le& \frac{B^2[(1+{b}/M)({M}/b)\lambda+C_0 d\lambda]} {
1-B[(1+{b}/M)({M}/b)\lambda+C_0 d\lambda]}.
\end{eqnarray*}
Since $d\lambda=o(1)$, we denote the right-hand side as $C_1 d\lambda$
for some $C_1>0$.

Next, we can further obtain a bound for $\|\hbtheta^{s.nd}_{(k)}-\btheta^*_{(k)}\|_{\ell_1}$.
\begin{eqnarray}\label{rankdantzigstep2}
\bigl\|\hbtheta^{s.nd}_{(k)}-\btheta^*_{(k)}
\bigr\|_{\ell_1} &\le& \bigl\| \bigl(\hat\theta_{kk}^{s.nd}-
\theta_{kk}^* \bigr)\hbeta^{s.nd}_{k}
\bigr\|_{\ell_1}+\bigl\| \theta_{kk}^* \bigl(\hbeta^{s.nd}_{k}-
\bbeta^*_{k} \bigr)\bigr\|_{\ell_1}
\nonumber
\\[-10pt]
\\[-10pt]
\nonumber
&\le& C_1 d\lambda\cdot b^{-1}M +B\cdot C_0
d\lambda.
\end{eqnarray}
Thus we can combine \eqref{rankdantzigstep1} and \eqref
{rankdantzigstep2} to derive the desired upper bound under the same
event. This completes the proof of Theorem~\ref{nonasymds1}.
\end{pf*}

\begin{pf*}{Proof of Theorem~\ref{nonasymads}}
Throughout the proof, we consider the event
\begin{equation}
\label{adsevent} \biggl\{\bigl\|\hbR^s-\bSigma^*\bigr\|_{\max}\le\min
\biggl(\lambda_{ad},\frac{b}M\lambda_{d} \biggr)
\biggr\}.
\end{equation}
For ease of notation, define $\lambda_{d}=\lambda_0$ and $\lambda_{ad}=\lambda_1$.
We focus on the proof of the sign consistency of $\hbeta^{s.nad}_k$ in
the sequel.

Under event \eqref{adsevent}, $\hbR^s_{\subaset\subaset}$ is always
positive-definite. To see this, the Weyl's inequality yields $\lambda_{\min}(\hbR^s_{\subaset\subaset})+\lambda_{\max}(\hbR^s_{\subaset
\subaset}-\bSigma^*_{\subaset\subaset})\ge\lambda_{\min}(\bSigma^*_{\subaset\subaset})$, and then we can bound the minimal eigenvalue
of $\hbR^s_{\subaset\subaset}$,
\[
\lambda_{\min} \bigl(\hbR^s_{\subaset\subaset} \bigr)\ge
\lambda_{\min} \bigl(\bSigma^*_{\subaset\subaset} \bigr)-\bigl\|
\hbR^s_{\subaset\subaset}-\bSigma^*_{\subaset
\subaset}\bigr\|_F \ge
\lambda_{1} \bigl(d-\sqrt{d(d-1)} \bigr)>0. 
\]
For each $k$ we introduce the dual variables $\balphaplus=(\alpha^+_j)_{j\neq k}\in\mathbb{R}^{p-1}_+$ and $\balphaminus=(\alpha^-_j)_{j\neq k}\in\mathbb{R}^{p-1}_+$. Then the Lagrange dual function
is defined as
\begin{eqnarray*}
L\bigl(\bbeta;\balphaplus,\balphaminus\bigr)&=&\bigl\|\bw^d_k\circ
\bbeta\bigr\|_{\ell
_1}+ \bigl(\hbR^s_{(k)}\bbeta-
\hbr^s_{(k)} -\lambda_{1}\bw^d_k
\bigr)^T\balphaplus
\\
&&{}+ \bigl(-\hbR^s_{(k)}\bbeta+\hbr^s_{(k)}-
\lambda_{1}\bw^d_k \bigr)^T
\balphaminus,
\end{eqnarray*}
where $\circ$ denotes the Hadamard product. Due to the strong duality
of linear programming [\citet{boyd2004}], the complementary slackness
condition holds for the primal problem with respect to any primal and
dual solution pair ($\bbeta,\balphaplus,\balphaminus$), which implies
that $
\alpha^+_j[(\hbR^s_{(k)}\bbeta-\hbr^s_{(k)})_j-\lambda_{1} w^d_j]=0
$
and
$
\alpha^-_j[-(\hbR^s_{(k)}\bbeta-\hbr^s_{(k)})_j-\lambda_{1} w^d_j]=0
$
for any $j\neq k$. Observe that only one of $\alpha^+_j$ and $\alpha^-_j$ can be zero since only one of
$(\hbR^s_{(k)}\bbeta-\hbr^s_{(k)})_j=\lambda_{1} w^d_j$ and $(\hbR^s_{(k)}\bbeta-\hbr^s_{(k)})_j=-\lambda_{1} w^d_j$ can hold indeed, and
thus we can uniquely define $\balpha_k=\balphaplus-\balphaminus$. Then
we can rewrite the Lagrange dual function as
\[
L(\bbeta;\balpha_k)= \bigl(\bw^d_k\circ
\sign(\bbeta)-\hbR^s_{(k)}\balpha_k
\bigr)^T\bbeta-\lambda_{1}\bigl\| \bw^d_k
\circ\balpha_k\bigr\|_{\ell_1}-\balpha_k^T
\hbr^s_{(k)}.
\]
By the Lagrange duality, the dual problem of \eqref{adadantzig} is
\[
\max_{\balpha\in\mathbb{R}^{p-1}}-\lambda_{1}\bigl\|\bw^d_k
\circ\balpha_k\bigr\|_{\ell_1}- \bigl\langle\balpha_k,
\hbr^s_{(k)} \bigr\rangle \qquad\mbox{subject to } \bigl|
\hbR^s_{(k)}\balpha_k\bigr|\le\bw^d_k.
\]
Now we shall construct an optimal primal and dual solution pair $(\tbeta_k,\talpha_k)$
to the rank-based adaptive Dantzig selector. In
addition, we show that $(\tbeta_k,\talpha_k)$ is actually the unique
solution pair to the rank-based adaptive Dantzig selector, and $\tbeta_k$ is
exactly supported in the true active set $\subaset$. To this
end, we construct $(\tbeta_k,\talpha_k)$ as
$
\talpha_k=(\talpha_{\subaset},\talpha_{\subaset^c})=(\talpha_{\subaset
},\bzero)
$
and
$
\tbeta_k=(\tbeta_{\subaset},\tbeta_{\subaset^c})=(\tbeta_{\subaset
},\bzero)
$
where $\talpha_{\subaset}=-(\hbR^s_{\subaset\subaset})^{-1}\bw^d_{\subaset}\circ\sign(\bbeta^*_{\subaset})$
and $\tbeta_{\subaset}=(\hbR^s_{\subaset\subaset})^{-1}(\hbr^s_{\subaset
}+\lambda_{1} \bw^d_{\subaset}\circ\sign(\talpha_{\subaset}))$.

In what follows, we first show that $(\tbeta_k,\talpha_k)$ satisfies
four optimality conditions, and then we will use these four optimality
conditions to prove that $(\tbeta_k,\talpha_k)$ is indeed a unique
optimal solution pair. The four optimality conditions are stated as follows:
\begin{eqnarray}
\hbR^s_{\subaset\subaset}\tbeta_{\subaset}-\hbr^s_{\subaset}
&=& \lambda_{1} \bw^d_{\subaset}\circ\sign(
\talpha_{\subaset}),\label{primalsupport}
\\
\hbR^s_{\subaset\subaset}\talpha_{\subaset}&=& -
\bw^d_{\subaset}\circ \sign(\tbeta_{\subaset}),
\label{dualsupport}
\\
\bigl|\hbR_{\subaset^c\subaset}^s\tbeta_{\subaset}-\hbr_{\subaset^c}^s\bigr|
&<& \lambda_{1} \bw^d_{\subaset^c},
\label{primalfeasible}
\\
\bigl|\hbR_{\subaset^c\subaset}^s\talpha_{\subaset}\bigr| &<&
\bw^d_{\subaset^c}, \label{dualfeasible}
\end{eqnarray}
where \eqref{primalsupport} and \eqref{primalfeasible} are primal
constraints, and \eqref{dualsupport} and \eqref{dualfeasible} are
dual constraints.

Note \eqref{primalsupport} can be easily verified by substituting
$\talpha_{\subaset}$ and $\tbeta_{\subaset}$.
Under \eqref{adsevent}, we can derive upper bounds for
$
K_1=\| (\hbR^s_{{\subaset\subaset}})^{-1}- (\bSigma^*_{{\subaset\subaset
}})^{-1}\|_{\ell_{\infty}}
$, and $
K_2=\| \hbR^s_{{\subaset^c\subaset}}(\hbR^s_{{\subaset\subaset}})^{-1}-
\bSigma^*_{{\subaset^c\subaset}}(\bSigma^*_{{\subaset\subaset}})^{-1}\|_{\ell_{\infty}}$.

Note $K_1= (\hbR^s_{{\subaset\subaset}})^{-1}\cdot(\hbR^s_{{\subaset
\subaset}}-\bSigma^*_{{\subaset\subaset}})\cdot(\bSigma^*_{{\subaset
\subaset}})^{-1}$, and then we have
\begin{eqnarray*}
K_1&\le&\bigl\| \bigl(\hbR^s_{{\subaset\subaset}}
\bigr)^{-1}\bigr\|_{\ell_{\infty}}\cdot \bigl\|\hbR^s_{{\subaset\subaset}}-
\bSigma^*_{{\subaset\subaset}}\bigr\|_{\ell
_{\infty}}\cdot\bigl\| \bigl(\bSigma^*_{{\subaset\subaset}}
\bigr)^{-1}\bigr\|_{\ell_{\infty
}}
\\
&\le& d\lambda_{1} G_k(G_k+K_1).
\end{eqnarray*}
Some simple calculation shows $
K_1\le\frac{d\lambda_{1} G_k^2}{1-d\lambda_{1} G_k}.
$
On the other hand,
\begin{eqnarray*}
K_2 &\le& \bigl\| \bigl(\hbR^s_{{\subaset^c\subaset}}-
\bSigma^*_{{\subaset^c\subaset}} \bigr) \bigl(\hbR^s_{{\subaset\subaset}}
\bigr)^{-1}\bigr\|_{\ell_{\infty}}
\\
&&{}+\bigl\|\bSigma^*_{{\subaset^c\subaset}} \bigl( \bigl(\bSigma^*_{{\subaset\subaset}}
\bigr)^{-1}- \bigl(\hbR^s_{{\subaset\subaset
}}
\bigr)^{-1} \bigr)\bigr\|_{\ell_{\infty}}
\\
&\le& \bigl(\bigl\|\hbR^s_{{\subaset^c\subaset}}-\bSigma^*_{{\subaset^c\subaset}}
\bigr\|_{\ell_{\infty}} +H_k \bigl\|\hbR^s_{{\subaset\subaset}}-
\bSigma^*_{{\subaset\subaset}}\bigr\|_{\ell
_{\infty}} \bigr)\cdot \bigl\| \bigl(
\hbR^s_{{\subaset\subaset}} \bigr)^{-1}\bigr\|_{\ell_{\infty}}
\\
&\le& d\lambda_{1}(1+H_k) (G_k+K_1)
\\
&\le& \frac{d\lambda_{1} G_k(1+H_k)}{1-d\lambda_{1} G_k}.
\end{eqnarray*}

Under probability event \eqref{adsevent}, we claim about $\bw^d_k$ that
\begin{eqnarray}
\bigl\|\bw^d_{\subaset^c}\bigr\|_{\min}&\ge& \frac{d\lambda_{1} G_k+H_k}{2\lambda_{1} G_k}
\psi_k+\frac{1+d G_k}{1-d\lambda_{1} G_k},\label{claim1}
\\
\bigl\|\bw^d_{\subaset}\bigr\|_{\infty}&\le& \frac{1-d\lambda_{1} G_k}{2\lambda_{1} G_k}
\psi_k-d G_k-1.\label{claim2}
\end{eqnarray}
This claim is very useful to prove the other three optimality
conditions \eqref{dualsupport}, \eqref{primalfeasible} and \eqref
{dualfeasible}, and their proofs will be provided later.

Now we are ready to prove \eqref{dualsupport}, \eqref{primalfeasible}
and \eqref{dualfeasible} for the solution pair $(\tbeta_k,\talpha_k)$.
To prove \eqref{dualsupport}, it is equivalent to show the sign
consistency that $\sign(\bbeta^*_{\subaset})
=
\sign(\tbeta_{\subaset})$ since we have $
\hbR^s_{\subaset\subaset}\talpha_{\subaset}
=
-\bw^d_{\subaset}\circ\sign(\bbeta^*_{\subaset})
$ if we plug in $\talpha_{\subaset}$ to its left-hand side. Recall that
$\bbeta^*_{\subaset}=(\bSigma^*_{{\subaset\subaset}})^{-1}\bsigma^*_{\subaset}$, and then we consider the difference between $\tbeta_{\subaset}$ and $\bbeta^*_{\subaset}$,
\begin{eqnarray*}
\tbeta_{\subaset}-\bbeta^*_{\subaset} &=& \bigl(\hbR^s_{\subaset\subaset}
\bigr)^{-1} \bigl(\hbr^s_{\subaset}-
\bsigma^*_{\subaset
}+ \lambda_{1} \bw^d_{\subaset}
\circ\sign( \talpha_{\subaset}) \bigr)
\\
&&{}- \bigl( \bigl(\hbR^s_{\subaset\subaset} \bigr)^{-1}-
\bigl( \bSigma^*_{{\subaset\subaset
}} \bigr)^{-1} \bigr)\bsigma^*_{\subaset}.
\end{eqnarray*}
Then we apply the triangle inequality to obtain an upper bound,
\begin{eqnarray*}
\bigl\|\tbeta_{\subaset}-\bbeta^*_{\subaset}\bigr\|_{\ell_\infty} &\le&
(G_k+K_1) \bigl(\lambda_{1}+
\lambda_{1}\bigl\|\bw^d_{\subaset}\bigr\|_{\ell_\infty}
\bigr)+K_1\bigl\|\bsigma^*_{\subaset}\bigr\|_{\ell_\infty}
\\
&\le& \frac{\lambda_{1} G_k}{1-d\lambda_{1} G_k} \bigl(1+\bigl\|\bw^d_{\subaset}
\bigr\|_{\ell
_\infty} \bigr)+\frac{d\lambda_{1} G_k^2}{1-d\lambda_{1} G_k}
\\
&<& \bigl\|\bbeta^*_{\subaset}\bigr\|_{\min},
\end{eqnarray*}
where the last inequality obviously holds by claim \eqref{claim1}. Then
by the above upper bound, $\sign(\bbeta^*_{\subaset})
=
\sign(\tbeta_{\subaset})$ will be immediately satisfied.

Next, we can easily obtain \eqref{dualfeasible} via the triangular inequality
\begin{eqnarray*}
\bigl\|\hbR_{\subaset^c\subaset}^s\talpha_{\subaset}\bigr\|_{\ell_\infty} &
\le& \bigl\|\hbR_{\subaset^c\subaset}^s \bigl(\hbR^s_{{\subaset\subaset}}
\bigr)^{-1}\bw^d_{\subaset}\bigr\|_{\ell_\infty}
\\
&\le& (H_k+K_2)\bigl\|\bw^d_{\subaset}
\bigr\|_{\infty}
\\
&\le& \frac{d\lambda_{1} G_k+H_k}{1-d\lambda_{1} G_k}\bigl\|\bw^d_{\subaset}\bigr\|_{\infty}
\\
&<& \bigl\|\bw^d_{\subaset^c}\bigr\|_{\min},
\end{eqnarray*}
where the last inequality can be easily shown by combining \eqref
{claim1} and \eqref{claim2}.

Now it remains to prove \eqref{primalfeasible}. Using the facts that
$\btheta^*_{{\subaset^c}}=\bzero$ and $\bSigma^*\bTheta^*=\bI$, simple
calculation yields that
$
\bSigma^*_{{\subaset^c\subaset}}(\bSigma^*_{{\subaset\subaset
}})^{-1}\bsigma^*_{\subaset}=\bsigma^*_{\subaset^c}.
$ Then we can rewrite the left-hand side of \eqref{primalfeasible} as
\begin{eqnarray*}
&&\hbR_{\subaset^c\subaset}^s\tbeta_{\subaset}-\hbr_{\subaset^c}^s
\\
&&\qquad=\hbR^s_{{\subaset^c\subaset}} \bigl(\hbR^s_{{\subaset\subaset}}
\bigr)^{-1} \bigl(\hbr^s_{\subaset}+
\lambda_{1} \bw^d_{\subaset}\circ\sign(
\talpha_{\subaset}) \bigr)-\hbr^s_{{\subaset^c}}
\\
&&\qquad=\hbR^s_{{\subaset^c\subaset}} \bigl(\hbR^s_{{\subaset\subaset}}
\bigr)^{-1} \bigl(\hbr^s_{\subaset}-
\bsigma^*_{\subaset
}+ \lambda_{1} \bw^d_{\subaset}
\circ\sign( \talpha_{\subaset}) \bigr)
\\
&&\qquad\quad{}+ \bigl(\hbR^s_{{\subaset^c\subaset}} \bigl(\hbR^s_{{\subaset\subaset}}
\bigr)^{-1} -\bSigma^*_{{\subaset^c\subaset}} \bigl(\bSigma^*_{{\subaset\subaset
}}
\bigr)^{-1} \bigr)\bsigma^*_{\subaset} + \bigl(\bsigma^*_{\subaset^c}-
\hbr_{\subaset^c}^s \bigr).
\end{eqnarray*}
Again we apply the triangle inequality to obtain an upper bound as follows:
\begin{eqnarray*}
&&\bigl\|\hbR_{\subaset^c\subaset}^s\tbeta_{\subaset}-
\hbr_{\subaset^c}^s\bigr\|_{\infty}
\\
&&\qquad\le (H_k+K_2) \bigl(\lambda_{1}+
\lambda_{1}\bigl\|\bw^d_{\subaset}\bigr\|_{\infty}
\bigr)+K_2\bigl\|\bsigma^*_{\subaset}\bigr\|_{\infty}+
\lambda_{1}
\\
&&\qquad\le \frac{d\lambda_{1}^2 G_k+\lambda_{1} H_k}{1-d\lambda_{1} G_k} \bigl(1+\bigl\| \bw^d_{\subaset}
\bigr\|_{\infty} \bigr) +\frac{d\lambda_{1} G_k(1+H_k)}{1-d\lambda_{1} G_k}+\lambda_{1}
\\
&&\qquad< \lambda_{1}\bigl\|\bw^d_{\subaset^c}\bigr\|_{\min},
\end{eqnarray*}
where the last inequality is due to \eqref{claim1} and \eqref{claim2}.


So far, the four optimality conditions have been verified for $(\tbeta_k,\talpha_k)$. In the sequel, we shall show that $\tbeta_k$ is indeed
a unique optimal solution. First, due to \eqref{primalfeasible} and
\eqref{dualfeasible}, $(\tbeta_k,\talpha_k)$ are feasible solutions to the
primal and dual problems, respectively. Then \eqref{primalsupport} and
\eqref{dualsupport} show that $(\tbeta_k,\talpha_k)$ satisfy the
complementary-slackness conditions for both the primal and the dual
problems. Thus, $(\tbeta_k,\talpha_k)$ are optimal solutions to these
problems by Theorem 4.5 in \citet{bertsimas1997}. Now it remains to show
the uniqueness. Suppose there exists another optimal solution $\brbeta_k$, and we have $\|\bw^d_k\circ\brbeta_k\|_{\ell_1}=\|\bw^d_k\circ
\tbeta_k\|_{\ell_1}$. Let $\Gamma_k$ denote the support of $\brbeta_k$,
and then $\brbeta_k=(\brbeta_{\Gamma_k},\bzero)$. By the strong duality
we have
\begin{eqnarray*}
\bigl\|\bw^d_k\circ\brbeta_k\bigr\|_{\ell_1}
&=&\bigl\|\bw^d_k\circ\tbeta_k\bigr\|_{\ell_1}
\\
&=&-\lambda_{1}\bigl\|\bw^d_k\circ
\talpha_k\bigr\|_{\ell_1}- \bigl\langle\talpha_k,
\hbr^s_{(k)} \bigr\rangle
\\
&=&\inf_{\bbeta}L\bigl(\bbeta;\talphaplus,\talphaminus\bigr)
\\
&\le& L\bigl(\brbeta_k;\talphaplus,\talphaminus\bigr)
\\
&\le& \bigl\|\bw^d_k\circ\brbeta_k
\bigr\|_{\ell_1}.
\end{eqnarray*}
Thus $L(\brbeta_k;\talphaplus,\talphaminus)=\|\bw^d_k\circ\brbeta_k\|_{\ell_1}$,
which immediately implies that the complementary slackness
condition holds for the primal problem, that is, $(\hbR^s_{(k)}\brbeta_k-\hbr^s_{(k)}-\lambda_{1}\bw^d_k)^T\talphaplus=0$
and $(-\hbR^s_{(k)}\brbeta_k+\hbr^s_{(k)}-\lambda_{1}\bw^d_k)^T\talphaminus=0$.
Now let $\brbetaplus=\max(\brbeta_k,\bzero)$ and $\brbetaminus=\min
(\brbeta_k,\bzero)$. Besides, we can similarly show that the
complementary slackness condition also holds for the dual problem, that
is, $(\hbR^s_{(k)}\talpha_k-\bw^d_k)^T\brbetaplus=0$ and $(-\hbR^s_{(k)}\talpha_k-\bw^d_k)^T\brbetaminus=0$.
Notice that $\talpha_{\subaset}\neq\bzero$ and $\talpha_{\subaset^c}=\bzero$ by
definition, and then we have
\begin{eqnarray}
\hbR^s_{\subaset\Gamma_k}\brbeta_{\Gamma_k}-\hbr^s_{\subaset}&=&
\lambda_{1} \bw^d_{\subaset}\circ\sign(
\talpha_{\subaset}),\label{unique5}
\\
\hbR^s_{\Gamma_k\subaset}\talpha_{\subaset}&=&-
\bw^d_{\Gamma_k}\circ \sign(\brbeta_{\Gamma_k}).
\label{unique6}
\end{eqnarray}
Observe that for any $j\in\Gamma_k$ but $j\notin\subaset$, $\hbR^s_{j\subaset}\talpha_{\subaset}=-w^d_j\sign(\breve\beta_j)$ in \eqref
{unique6} cannot hold since it contradicts with \eqref{dualfeasible}.
Then it is easy to see that $\Gamma_k\subset\subaset$ obviously holds,
which immediately implies that $\hbeta_{\subaset}$ and $\brbeta_{\Gamma
_k}$ satisfy the same optimality condition \eqref{primalsupport}. Thus
the uniqueness follows from \eqref{primalsupport}, \eqref{unique5} and
the nonsingularity of $\hbR^s_{\subaset\subaset}$.

Now it remains to verify the claims \eqref{claim1} and \eqref{claim2}
under event \eqref{adsevent}. Under the event $\|\hbR^s-\bSigma^*\|_{\max}\le b\lambda_0/M$, it has been shown in Theorem~\ref{nonasymds1} that for some
$C_0=4B^2(2+\frac{b}M)>0$, we have $\|\hbeta^{s.nd}_{k}-\bbeta^*_k\|_{\ell_1}\le C_0  d\lambda_0$. Then we can derive a lower bound for $\|
\bw^d_{\subaset^c}\|_{\min}$,
\[
\bigl\|\bw^d_{\subaset^c}\bigr\|_{\min}=\frac{1}{\max_{j\in\subaset^c}|\hat
\beta^{s.nd}_j|+{1}/n} \ge
\frac{1}{C_0 d\lambda_0+{1}/n},
\]
which immediately yields the desired lower bound by noting that
\[
\frac{G_kd\lambda_1 +H_k}{2G_k\cdot\lambda_1 }\cdot\psi_{k}+\frac{
1+G_k d}{1- G_k d\lambda_1} \le
\frac{H_k\psi_{k}}{2G_k\cdot\lambda_1}+(\psi_k+2G_k)\cdot d+2 \le
\frac{1}{C_0 d\lambda_0+{1}/n},
\]
where both inequalities follow from the proper choices of tuning
parameters $\lambda_0$ and $\lambda_1$ as stated in Theorem~\ref
{nonasymads}. On the other hand,
\[
\frac{1-G_k\cdot d\lambda_1 }{2G_k\cdot\lambda_1 }\psi_{k}-d G_k-1 \ge
\frac{\psi_{k}}{2G_k\cdot\lambda_1 }-(\psi_k+G_k)\cdot d-1 \ge
\frac{\psi_{k}}{4G_k\cdot\lambda_1},
\]
where the last inequality follows from the proper choice of $\lambda_1$
as stated in Theorem~\ref{nonasymads}. Likewise we can prove the
second claim (\ref{claim2}) by noticing that
\[
\bigl\|\bw^d_{\subaset}\bigr\|_{\infty} \le\frac{1}{\min_{j\in\subaset}|\hat\beta^{s.nd}_j|} \le
\frac{1}{\psi_{k}-C_0 d\lambda_0} \le\frac{2}{\psi_k} \le\frac{\psi_{k}}{4G_k\cdot\lambda_1},
\]
where we use facts that $\psi_{k}\ge2C_0 d\lambda_0$ and $\psi^2_{k}\ge
8G_k\lambda_1$. The two claims are proved, which completes the proof of
Theorem~\ref{nonasymads}.\vadjust{\goodbreak}
\end{pf*}


\begin{pf*}{Proof of Theorem~\ref{nonasymclime}}
To bound the difference between $\hbTheta_c^s$ and $\bTheta^*$ under
the entry-wise $\ell_\infty$-norm, we consider the event $\{\|\hbR^s-\bSigma\|_{\max}\le\lambda/M\}$. First, we show that $\bTheta^*$ is
always a feasible solution under the above event,
\[
\bigl\|\hbR^s\bTheta^*-\bI\bigr\|_{\max} \le \bigl\| \bigl(
\hbR^s-\bSigma^* \bigr)\bTheta^*\bigr\|_{\max} \le \bigl\|
\hbR^s-\bSigma\bigr\|_{\max}\cdot\bigl\|\bTheta^*\bigr\|_{\ell_1} \le
\lambda.
\]
Note that $\hbTheta_c^s$ is the optimal solution, and then $\|\hbR^s\hbTheta_c^s-\bI\|_{\max}\le\lambda$ obviously holds. Moreover, it
is easy to see that by definition $\|\hbTheta_c^s\|_{\ell_1}\le\|\bTheta^*\|_{\ell_1}$ always holds. Now we can obtain the desired bound under
the entry-wise $\ell_\infty$-norm.
\begin{eqnarray*}
\bigl\|\hbTheta_c^s-\bTheta^*\bigr\|_{\max} &\le& \bigl\|
\bTheta^*\bigr\|_{\ell_1}\cdot\bigl\|\bSigma^*\hbTheta_c^s-\bI
\bigr\|_{\max}
\\
&=& M\cdot\bigl\| \bigl(\bSigma^*-\hbR^s \bigr)\hbTheta^s+
\hbR^s\hbTheta_c^s-\bI\bigr\|_{\max}
\\
&\le& M\cdot\bigl\|\bSigma^*-\hbR^s\bigr\|_{\max}\cdot\bigl\|
\hbTheta^s_c\bigr\|_{\ell_1}+M\cdot\bigl\| \hbR^s
\hbTheta_c^s-\bI\bigr\|_{\max}
\\
&\le& \lambda\cdot\bigl\|\bTheta^*\bigr\|_{\ell_1}+M\lambda
\\
&=& 2M\lambda.
\end{eqnarray*}
\upqed\end{pf*}

\begin{pf*}{Proof of Theorem~\ref{nonasymaclime}}
The techniques we use are similar to these for the proof of Theorem \ref
{nonasymads}. The detailed proof of Theorem~\ref{nonasymaclime} is
relegated to the supplementary material [\citet{supp}] and also the technical report
version of this paper [\citet{XZ11b}] for the sake of space.
\end{pf*}
\end{appendix}


\section*{\texorpdfstring{Acknowledgments.}{Acknowledgments}}
We thank the Editor, the Associate Editor and three referees for their
helpful comments.

\begin{supplement}
\stitle{Supplement material for ``Regularized rank-based estimation of
high-dimensional nonparanormal graphical models''\\}
\slink[doi,text={10.1214/12-\break AOS1041SUPP}]{10.1214/12-AOS1041SUPP} 
\sdatatype{.pdf}
\sfilename{aos1041\_supp.pdf}
\sdescription{In this supplementary note, we give the complete proofs
of Theorems~\ref{nonasymds1} and~\ref{nonasymaclime}.}
\end{supplement}


%
%

\printaddresses


\begin{thebibliography}{45}

\bibitem[\protect\citeauthoryear{Banerjee, El~Ghaoui and
d'Aspremont}{2008}]{banerjee2008}
\begin{barticle}[mr]
\bauthor{\bsnm{Banerjee},~\bfnm{Onureena}\binits{O.}},
\bauthor{\bsnm{El~Ghaoui},~\bfnm{Laurent}\binits{L.}} \AND
\bauthor{\bsnm{d'Aspremont},~\bfnm{Alexandre}\binits{A.}}
(\byear{2008}).
\btitle{Model selection through sparse maximum likelihood estimation for
multivariate {G}aussian or binary data}.
\bjournal{J.~Mach. Learn. Res.}
\bvolume{9}
\bpages{485--516}.
\bid{issn={1532-4435}, mr={2417243}}
\bptok{imsref}%
\end{barticle}
\endbibitem

\bibitem[\protect\citeauthoryear{Bertsimas and
Tsitsiklis}{1997}]{bertsimas1997}
\begin{bbook}[author]
\bauthor{\bsnm{Bertsimas},~\bfnm{D.}\binits{D.}} \AND
\bauthor{\bsnm{Tsitsiklis},~\bfnm{J.~N.}\binits{J.~N.}}
(\byear{1997}).
\btitle{Introduction to Linear Optimization}.
\bpublisher{Athena Scientific}, \blocation{Belmont, MA}.
\bptok{imsref}%
\end{bbook}
\endbibitem

\bibitem[\protect\citeauthoryear{Bickel, Ritov and Tsybakov}{2009}]{dasso1}
\begin{barticle}[mr]
\bauthor{\bsnm{Bickel},~\bfnm{Peter~J.}\binits{P.~J.}},
\bauthor{\bsnm{Ritov},~\bfnm{Ya'acov}\binits{Y.}} \AND
\bauthor{\bsnm{Tsybakov},~\bfnm{Alexandre~B.}\binits{A.~B.}}
(\byear{2009}).
\btitle{Simultaneous analysis of lasso and {D}antzig selector}.
\bjournal{Ann. Statist.}
\bvolume{37}
\bpages{1705--1732}.
\bid{doi={10.1214/08-AOS620}, issn={0090-5364}, mr={2533469}}
\bptok{imsref}%
\end{barticle}
\endbibitem

\bibitem[\protect\citeauthoryear{Boyd and Vandenberghe}{2004}]{boyd2004}
\begin{bbook}[mr]
\bauthor{\bsnm{Boyd},~\bfnm{Stephen}\binits{S.}} \AND
\bauthor{\bsnm{Vandenberghe},~\bfnm{Lieven}\binits{L.}}
(\byear{2004}).
\btitle{Convex Optimization}.
\bpublisher{Cambridge Univ. Press}, \blocation{Cambridge}.
\bid{mr={2061575}}
\bptok{imsref}%
\end{bbook}
\endbibitem

\bibitem[\protect\citeauthoryear{Cai, Liu and Luo}{2011}]{clime11}
\begin{barticle}[mr]
\bauthor{\bsnm{Cai},~\bfnm{Tony}\binits{T.}},
\bauthor{\bsnm{Liu},~\bfnm{Weidong}\binits{W.}} \AND
\bauthor{\bsnm{Luo},~\bfnm{Xi}\binits{X.}}
(\byear{2011}).
\btitle{A constrained {$\ell\sb 1$} minimization approach to sparse precision
matrix estimation}.
\bjournal{J. Amer. Statist. Assoc.}
\bvolume{106}
\bpages{594--607}.
\bid{doi={10.1198/jasa.2011.tm10155}, issn={0162-1459}, mr={2847973}}
\bptok{imsref}%
\end{barticle}
\endbibitem

\bibitem[\protect\citeauthoryear{Candes and Tao}{2007}]{candes2007}
\begin{barticle}[mr]
\bauthor{\bsnm{Candes},~\bfnm{Emmanuel}\binits{E.}} \AND
\bauthor{\bsnm{Tao},~\bfnm{Terence}\binits{T.}}
(\byear{2007}).
\btitle{The {D}antzig selector: Statistical estimation when {$p$} is much
larger than {$n$}}.
\bjournal{Ann. Statist.}
\bvolume{35}
\bpages{2313--2351}.
\bid{doi={10.1214/009053606000001523}, issn={0090-5364}, mr={2382644}}
\bptok{imsref}%
\end{barticle}\vadjust{\goodbreak}
\endbibitem

\bibitem[\protect\citeauthoryear{Chen and Fan}{2006}]{chen2006}
\begin{barticle}[mr]
\bauthor{\bsnm{Chen},~\bfnm{Xiaohong}\binits{X.}} \AND
\bauthor{\bsnm{Fan},~\bfnm{Yanqin}\binits{Y.}}
(\byear{2006}).
\btitle{Estimation of copula-based semiparametric time series models}.
\bjournal{J.~Econometrics}
\bvolume{130}
\bpages{307--335}.
\bid{doi={10.1016/j.jeconom.2005.03.004}, issn={0304-4076}, mr={2211797}}
\bptok{imsref}%
\end{barticle}
\endbibitem

\bibitem[\protect\citeauthoryear{Chen, Fan and
Tsyrennikov}{2006}]{chen2006jasa}
\begin{barticle}[mr]
\bauthor{\bsnm{Chen},~\bfnm{Xiaohong}\binits{X.}},
\bauthor{\bsnm{Fan},~\bfnm{Yanqin}\binits{Y.}} \AND
\bauthor{\bsnm{Tsyrennikov},~\bfnm{Viktor}\binits{V.}}
(\byear{2006}).
\btitle{Efficient estimation of semiparametric multivariate copula models}.
\bjournal{J. Amer. Statist. Assoc.}
\bvolume{101}
\bpages{1228--1240}.
\bid{doi={10.1198/016214506000000311}, issn={0162-1459}, mr={2328309}}
\bptok{imsref}%
\end{barticle}
\endbibitem

\bibitem[\protect\citeauthoryear{Dempster}{1972}]{dempster1972}
\begin{barticle}[author]
\bauthor{\bsnm{Dempster},~\bfnm{AP}\binits{A.}}
(\byear{1972}).
\btitle{{Covariance selection}}.
\bjournal{Biometrics}
\bvolume{28}
\bpages{157--175}.
\bptok{imsref}%
\end{barticle}
\endbibitem

\bibitem[\protect\citeauthoryear{Devlin, Gnanadesikan and
Kettenring}{1975}]{devlin1975}
\begin{barticle}[author]
\bauthor{\bsnm{Devlin},~\bfnm{S.~J.}\binits{S.~J.}},
\bauthor{\bsnm{Gnanadesikan},~\bfnm{R.}\binits{R.}} \AND
\bauthor{\bsnm{Kettenring},~\bfnm{J.~R.}\binits{J.~R.}}
(\byear{1975}).
\btitle{Robust estimation and outlier detection with correlation coefficients}.
\bjournal{Biometrika}
\bvolume{62}
\bpages{531--545}.
\bptok{imsref}%
\end{barticle}
\endbibitem

\bibitem[\protect\citeauthoryear{Dicker and Lin}{2009}]{dicker2009}
\begin{bmisc}[author]
\bauthor{\bsnm{Dicker},~\bfnm{L.}\binits{L.}} \AND
\bauthor{\bsnm{Lin},~\bfnm{X.}\binits{X.}}
(\byear{2009}).
\bhowpublished{Variable selection using the Dantzig selector: Asymptotic theory and
extensions. Unpublished manuscript.}
\bptok{imsref}%
\end{bmisc}
\endbibitem

\bibitem[\protect\citeauthoryear{Dobra, Eicher and Lenkoski}{2010}]{dobra2009}
\begin{barticle}[mr]
\bauthor{\bsnm{Dobra},~\bfnm{Adrian}\binits{A.}},
\bauthor{\bsnm{Eicher},~\bfnm{Theo~S.}\binits{T.~S.}} \AND
\bauthor{\bsnm{Lenkoski},~\bfnm{Alex}\binits{A.}}
(\byear{2010}).
\btitle{Modeling uncertainty in macroeconomic growth determinants using
{G}aussian graphical models}.
\bjournal{Stat. Methodol.}
\bvolume{7}
\bpages{292--306}.
\bid{doi={10.1016/j.stamet.2009.11.003}, issn={1572-3127}, mr={2643603}}
\bptnote{check year}%
\bptok{imsref}%
\end{barticle}
\endbibitem

\bibitem[\protect\citeauthoryear{Drton and Perlman}{2004}]{drton2004}
\begin{barticle}[mr]
\bauthor{\bsnm{Drton},~\bfnm{Mathias}\binits{M.}} \AND
\bauthor{\bsnm{Perlman},~\bfnm{Michael~D.}\binits{M.~D.}}
(\byear{2004}).
\btitle{Model selection for {G}aussian concentration graphs}.
\bjournal{Biometrika}
\bvolume{91}
\bpages{591--602}.
\bid{doi={10.1093/biomet/91.3.591}, issn={0006-3444}, mr={2090624}}
\bptok{imsref}%
\end{barticle}
\endbibitem

\bibitem[\protect\citeauthoryear{Drton and Perlman}{2007}]{drton2007}
\begin{barticle}[mr]
\bauthor{\bsnm{Drton},~\bfnm{Mathias}\binits{M.}} \AND
\bauthor{\bsnm{Perlman},~\bfnm{Michael~D.}\binits{M.~D.}}
(\byear{2007}).
\btitle{Multiple testing and error control in {G}aussian graphical model
selection}.
\bjournal{Statist. Sci.}
\bvolume{22}
\bpages{430--449}.
\bid{doi={10.1214/088342307000000113}, issn={0883-4237}, mr={2416818}}
\bptok{imsref}%
\end{barticle}
\endbibitem

\bibitem[\protect\citeauthoryear{Edwards}{2000}]{edwards2000}
\begin{bbook}[mr]
\bauthor{\bsnm{Edwards},~\bfnm{David}\binits{D.}}
(\byear{2000}).
\btitle{Introduction to Graphical Modelling},
\bedition{2nd} ed.
\bpublisher{Springer}, \blocation{New York}.
\bid{doi={10.1007/978-1-4612-0493-0}, mr={1880319}}
\bptok{imsref}%
\end{bbook}
\endbibitem

\bibitem[\protect\citeauthoryear{Fan and Li}{2001}]{fan2001}
\begin{barticle}[mr]
\bauthor{\bsnm{Fan},~\bfnm{Jianqing}\binits{J.}} \AND
\bauthor{\bsnm{Li},~\bfnm{Runze}\binits{R.}}
(\byear{2001}).
\btitle{Variable selection via nonconcave penalized likelihood and its oracle
properties}.
\bjournal{J. Amer. Statist. Assoc.}
\bvolume{96}
\bpages{1348--1360}.
\bid{doi={10.1198/016214501753382273}, issn={0162-1459}, mr={1946581}}
\bptok{imsref}%
\end{barticle}
\endbibitem

\bibitem[\protect\citeauthoryear{Friedman}{2004}]{nfriedman2004}
\begin{barticle}[pbm]
\bauthor{\bsnm{Friedman},~\bfnm{Nir}\binits{N.}}
(\byear{2004}).
\btitle{Inferring cellular networks using probabilistic graphical models}.
\bjournal{Science}
\bvolume{303}
\bpages{799--805}.
\bid{doi={10.1126/science.1094068}, issn={1095-9203}, pii={303/5659/799},
pmid={14764868}}
\bptok{imsref}%
\end{barticle}
\endbibitem

\bibitem[\protect\citeauthoryear{Friedman, Hastie and
Tibshirani}{2008}]{friedman2008}
\begin{barticle}[pbm]
\bauthor{\bsnm{Friedman},~\bfnm{Jerome}\binits{J.}},
\bauthor{\bsnm{Hastie},~\bfnm{Trevor}\binits{T.}} \AND
\bauthor{\bsnm{Tibshirani},~\bfnm{Robert}\binits{R.}}
(\byear{2008}).
\btitle{Sparse inverse covariance estimation with the graphical lasso}.
\bjournal{Biostatistics}
\bvolume{9}
\bpages{432--441}.
\bid{doi={10.1093/biostatistics/kxm045}, issn={1468-4357}, mid={NIHMS248717},
pii={kxm045}, pmcid={3019769}, pmid={18079126}}
\bptok{imsref}%
\end{barticle}
\endbibitem

\bibitem[\protect\citeauthoryear{Hoeffding}{1948}]{hoeffding1948}
\begin{barticle}[mr]
\bauthor{\bsnm{Hoeffding},~\bfnm{Wassily}\binits{W.}}
(\byear{1948}).
\btitle{A class of statistics with asymptotically normal distribution}.
\bjournal{Ann. Math. Statistics}
\bvolume{19}
\bpages{293--325}.
\bid{issn={0003-4851}, mr={0026294}}
\bptok{imsref}%
\end{barticle}
\endbibitem

\bibitem[\protect\citeauthoryear{James, Radchenko and Lv}{2009}]{dasso2009}
\begin{barticle}[mr]
\bauthor{\bsnm{James},~\bfnm{Gareth~M.}\binits{G.~M.}},
\bauthor{\bsnm{Radchenko},~\bfnm{Peter}\binits{P.}} \AND
\bauthor{\bsnm{Lv},~\bfnm{Jinchi}\binits{J.}}
(\byear{2009}).
\btitle{D{ASSO}: Connections between the {D}antzig selector and lasso}.
\bjournal{J. R. Stat. Soc. Ser. B Stat. Methodol.}
\bvolume{71}
\bpages{127--142}.
\bid{doi={10.1111/j.1467-9868.2008.00668.x}, issn={1369-7412}, mr={2655526}}
\bptok{imsref}%
\end{barticle}
\endbibitem

\bibitem[\protect\citeauthoryear{Kendall}{1948}]{kendall1948}
\begin{bbook}[author]
\bauthor{\bsnm{Kendall},~\bfnm{M.~G.}\binits{M.~G.}}
(\byear{1948}).
\btitle{Rank Correlation Methods}.
\bpublisher{Charles Griffin and Co. Ltd.}, \blocation{London}.
\bptok{imsref}%
\end{bbook}
\endbibitem

\bibitem[\protect\citeauthoryear{Kruskal}{1958}]{kruskal1958}
\begin{barticle}[mr]
\bauthor{\bsnm{Kruskal},~\bfnm{William~H.}\binits{W.~H.}}
(\byear{1958}).
\btitle{Ordinal measures of association}.
\bjournal{J. Amer. Statist. Assoc.}
\bvolume{53}
\bpages{814--861}.
\bid{issn={0162-1459}, mr={0100941}}
\bptok{imsref}%
\end{barticle}
\endbibitem

\bibitem[\protect\citeauthoryear{Lam and Fan}{2009}]{fan2009a}
\begin{barticle}[mr]
\bauthor{\bsnm{Lam},~\bfnm{Clifford}\binits{C.}} \AND
\bauthor{\bsnm{Fan},~\bfnm{Jianqing}\binits{J.}}
(\byear{2009}).
\btitle{Sparsistency and rates of convergence in large covariance matrix
estimation}.
\bjournal{Ann. Statist.}
\bvolume{37}
\bpages{4254--4278}.
\bid{doi={10.1214/09-AOS720}, issn={0090-5364}, mr={2572459}}
\bptok{imsref}%
\end{barticle}
\endbibitem

\bibitem[\protect\citeauthoryear{Laule et~al.}{2003}]{laule2003}
\begin{barticle}[author]
\bauthor{\bsnm{Laule},~\bfnm{O.}\binits{O.}},
\bauthor{\bsnm{F{\"u}rholz},~\bfnm{A.}\binits{A.}},
\bauthor{\bsnm{Chang},~\bfnm{H.~S.}\binits{H.~S.}},
\bauthor{\bsnm{Zhu},~\bfnm{T.}\binits{T.}},
\bauthor{\bsnm{Wang},~\bfnm{X.}\binits{X.}},
\bauthor{\bsnm{Heifetz},~\bfnm{P.~B.}\binits{P.~B.}},
\bauthor{\bsnm{Gruissem},~\bfnm{W.}\binits{W.}} \AND
\bauthor{\bsnm{Lange},~\bfnm{M.}\binits{M.}}
(\byear{2003}).
\btitle{{Crosstalk between cytosolic and plastidial pathways of isoprenoid
biosynthesis in Arabidopsis thaliana}}.
\bjournal{Proc. Natl. Acad. Sci. USA}
\bvolume{100}
\bpages{6866--6871}.
\bptok{imsref}%
\end{barticle}
\endbibitem

\bibitem[\protect\citeauthoryear{Lauritzen}{1996}]{lauritzen1996}
\begin{bbook}[mr]
\bauthor{\bsnm{Lauritzen},~\bfnm{Steffen~L.}\binits{S.~L.}}
(\byear{1996}).
\btitle{Graphical Models}.
\bseries{Oxford Statistical Science Series}
\bvolume{17}.
\bpublisher{The Clarendon Press Oxford Univ. Press}, \blocation{New York}.
\bid{mr={1419991}}
\bptok{imsref}%
\end{bbook}
\endbibitem

\bibitem[\protect\citeauthoryear{Lehmann}{1998}]{lehmann1998}
\begin{bbook}[author]
\bauthor{\bsnm{Lehmann},~\bfnm{E.~L.}\binits{E.~L.}}
(\byear{1998}).
\btitle{{Nonparametrics: Statistical Methods Based on Ranks}}.
\bpublisher{Prentice Hall Upper Saddle River}, \blocation{New Jersey}.
\bptok{imsref}%
\end{bbook}
\endbibitem

\bibitem[\protect\citeauthoryear{Li and Gui}{2006}]{lihz2006}
\begin{barticle}[pbm]
\bauthor{\bsnm{Li},~\bfnm{Hongzhe}\binits{H.}} \AND
\bauthor{\bsnm{Gui},~\bfnm{Jiang}\binits{J.}}
(\byear{2006}).
\btitle{Gradient directed regularization for sparse Gaussian concentration
graphs, with applications to inference of genetic networks}.
\bjournal{Biostatistics}
\bvolume{7}
\bpages{302--317}.
\bid{doi={10.1093/biostatistics/kxj008}, issn={1465-4644}, pii={kxj008},
pmid={16326758}}
\bptok{imsref}%
\end{barticle}
\endbibitem

\bibitem[\protect\citeauthoryear{Liu, Lafferty and
Wasserman}{2009}]{lafferty2009}
\begin{barticle}[mr]
\bauthor{\bsnm{Liu},~\bfnm{Han}\binits{H.}},
\bauthor{\bsnm{Lafferty},~\bfnm{John}\binits{J.}} \AND
\bauthor{\bsnm{Wasserman},~\bfnm{Larry}\binits{L.}}
(\byear{2009}).
\btitle{The nonparanormal: Semiparametric estimation of high dimensional
undirected graphs}.
\bjournal{J. Mach. Learn. Res.}
\bvolume{10}
\bpages{2295--2328}.
\bid{issn={1532-4435}, mr={2563983}}
\bptok{imsref}%
\end{barticle}
\endbibitem

\bibitem[\protect\citeauthoryear{Liu et~al.}{2012}]{liu2012}
\begin{bmisc}[author]
\bauthor{\bsnm{Liu},~\bfnm{H.}\binits{H.}},
\bauthor{\bsnm{Han},~\bfnm{F.}\binits{F.}},
\bauthor{\bsnm{Yuan},~\bfnm{M.}\binits{M.}},
\bauthor{\bsnm{Lafferty},~\bfnm{J.}\binits{J.}} \AND
\bauthor{\bsnm{Wasserman},~\bfnm{L.}\binits{L.}}
(\byear{2012}).
\bhowpublished{High dimensional semiparametric Gaussian copula graphical models. Technical report, Johns Hopkins Univ.}
\bptok{imsref}%
\end{bmisc}
\endbibitem

\bibitem[\protect\citeauthoryear{McDiarmid}{1989}]{mcdiarmid1989}
\begin{bincollection}[mr]
\bauthor{\bsnm{McDiarmid},~\bfnm{Colin}\binits{C.}}
(\byear{1989}).
\btitle{On the method of bounded differences}.
In \bbooktitle{Surveys in Combinatorics, 1989 ({N}orwich, 1989)}.
\bseries{London Mathematical Society Lecture Note Series}
\bvolume{141}
\bpages{148--188}.
\bpublisher{Cambridge Univ. Press}, \blocation{Cambridge}.
\bid{mr={1036755}}
\bptok{imsref}%
\end{bincollection}
\endbibitem

\bibitem[\protect\citeauthoryear{Meinshausen and
B{\"u}hlmann}{2006}]{meinshausen2006}
\begin{barticle}[mr]
\bauthor{\bsnm{Meinshausen},~\bfnm{Nicolai}\binits{N.}} \AND
\bauthor{\bsnm{B{\"u}hlmann},~\bfnm{Peter}\binits{P.}}
(\byear{2006}).
\btitle{High-dimensional graphs and variable selection with the lasso}.
\bjournal{Ann. Statist.}
\bvolume{34}
\bpages{1436--1462}.
\bid{doi={10.1214/009053606000000281}, issn={0090-5364}, mr={2278363}}
\bptok{imsref}%
\end{barticle}
\endbibitem

\bibitem[\protect\citeauthoryear{Peng et al.}{2009}]{peng2009}
\begin{barticle}[mr]
\bauthor{\bsnm{Peng},~\bfnm{Jie}\binits{J.}},
\bauthor{\bsnm{Wang},~\bfnm{Pei}\binits{P.}},
\bauthor{\bsnm{Zhou},~\bfnm{Nengfeng}\binits{N.}} \AND
\bauthor{\bsnm{Zhu},~\bfnm{Ji}\binits{J.}}
(\byear{2009}).
\btitle{Partial correlation estimation by joint sparse regression models}.
\bjournal{J. Amer. Statist. Assoc.}
\bvolume{104}
\bpages{735--746}.
\bid{doi={10.1198/jasa.2009.0126}, issn={0162-1459}, mr={2541591}}
\bptok{imsref}%
\end{barticle}
\endbibitem

\bibitem[\protect\citeauthoryear{Ravikumar et~al.}{2011}]{ravikumar2008}
\begin{barticle}[author]
\bauthor{\bsnm{Ravikumar},~\bfnm{P.}\binits{P.}},
\bauthor{\bsnm{Wainwright},~\bfnm{M.~J.}\binits{M.~J.}},
\bauthor{\bsnm{Raskutti},~\bfnm{G.}\binits{G.}} \AND
\bauthor{\bsnm{Yu},~\bfnm{B.}\binits{B.}}
(\byear{2011}).
\btitle{High-dimensional covariance estimation by minimizing $\ell_1$-penalized log-determinant
divergence}.
\bjournal{Electron. J. Statist.}
\bvolume{5}
\bpages{935--980}.
\bptok{imsref}%
\end{barticle}
\endbibitem

\bibitem[\protect\citeauthoryear{Rodr{\'{\i}}guez-Concepci{\'{o}}n
et~al.}{2004}]{rodriguez2004}
\begin{barticle}[pbm]
\bauthor{\bsnm{Rodr{\'{\i}}guez-Concepci{\'{o}}n},~\bfnm{Manuel}\binits{M.}},
\bauthor{\bsnm{For{\'{e}}s},~\bfnm{Oriol}\binits{O.}},
\bauthor{\bsnm{Martinez-Garc{\'{\i}}a},~\bfnm{Jaime~F.}\binits{J.~F.}},
\bauthor{\bsnm{Gonz{\'{a}}lez},~\bfnm{Victor}\binits{V.}},
\bauthor{\bsnm{Phillips},~\bfnm{Michael~A.}\binits{M.~A.}},
\bauthor{\bsnm{Ferrer},~\bfnm{Albert}\binits{A.}} \AND
\bauthor{\bsnm{Boronat},~\bfnm{Albert}\binits{A.}}
(\byear{2004}).
\btitle{Distinct light-mediated pathways regulate the biosynthesis and exchange
of isoprenoid precursors during Arabidopsis seedling development}.
\bjournal{Plant Cell}
\bvolume{16}
\bpages{144--156}.
\bid{doi={10.1105/tpc.016204}, issn={1040-4651}, pii={tpc.016204},
pmcid={301401}, pmid={14660801}}
\bptok{imsref}%
\end{barticle}
\endbibitem

\bibitem[\protect\citeauthoryear{Rothman et~al.}{2008}]{rothman2008}
\begin{barticle}[mr]
\bauthor{\bsnm{Rothman},~\bfnm{Adam~J.}\binits{A.~J.}},
\bauthor{\bsnm{Bickel},~\bfnm{Peter~J.}\binits{P.~J.}},
\bauthor{\bsnm{Levina},~\bfnm{Elizaveta}\binits{E.}} \AND
\bauthor{\bsnm{Zhu},~\bfnm{Ji}\binits{J.}}
(\byear{2008}).
\btitle{Sparse permutation invariant covariance estimation}.
\bjournal{Electron. J. Stat.}
\bvolume{2}
\bpages{494--515}.
\bid{doi={10.1214/08-EJS176}, issn={1935-7524}, mr={2417391}}
\bptok{imsref}%
\end{barticle}
\endbibitem

\bibitem[\protect\citeauthoryear{Song}{2000}]{song2000}
\begin{barticle}[mr]
\bauthor{\bsnm{Song},~\bfnm{Peter Xue-Kun}\binits{P.~X.-K.}}
(\byear{2000}).
\btitle{Multivariate dispersion models generated from {G}aussian copula}.
\bjournal{Scand. J. Stat.}
\bvolume{27}
\bpages{305--320}.
\bid{doi={10.1111/1467-9469.00191}, issn={0303-6898}, mr={1777506}}
\bptok{imsref}%
\end{barticle}
\endbibitem

\bibitem[\protect\citeauthoryear{Tibshirani}{1996}]{tibshirani1996}
\begin{barticle}[mr]
\bauthor{\bsnm{Tibshirani},~\bfnm{Robert}\binits{R.}}
(\byear{1996}).
\btitle{Regression shrinkage and selection via the lasso}.
\bjournal{J. Roy. Statist. Soc. Ser. B}
\bvolume{58}
\bpages{267--288}.
\bid{issn={0035-9246}, mr={1379242}}
\bptok{imsref}%
\end{barticle}
\endbibitem

\bibitem[\protect\citeauthoryear{van~de Geer, B{\"u}hlmann and
Zhou}{2011}]{Geer2011}
\begin{barticle}[mr]
\bauthor{\bparticle{van~de} \bsnm{Geer},~\bfnm{Sara}\binits{S.}},
\bauthor{\bsnm{B{\"u}hlmann},~\bfnm{Peter}\binits{P.}} \AND
\bauthor{\bsnm{Zhou},~\bfnm{Shuheng}\binits{S.}}
(\byear{2011}).
\btitle{The adaptive and the thresholded {L}asso for potentially misspecified
models (and a lower bound for the {L}asso)}.
\bjournal{Electron. J. Stat.}
\bvolume{5}
\bpages{688--749}.
\bid{doi={10.1214/11-EJS624}, issn={1935-7524}, mr={2820636}}
\bptok{imsref}%
\end{barticle}
\endbibitem

\bibitem[\protect\citeauthoryear{Wille et~al.}{2004}]{wille2004}
\begin{barticle}[author]
\bauthor{\bsnm{Wille},~\bfnm{A.}\binits{A.}},
\bauthor{\bsnm{Zimmermann},~\bfnm{P.}\binits{P.}},
\bauthor{\bsnm{Vranov{\'a}},~\bfnm{E.}\binits{E.}},
\bauthor{\bsnm{F{\"u}rholz},~\bfnm{A.}\binits{A.}},
\bauthor{\bsnm{Laule},~\bfnm{O.}\binits{O.}},
\bauthor{\bsnm{Bleuler},~\bfnm{S.}\binits{S.}},
\bauthor{\bsnm{Hennig},~\bfnm{L.}\binits{L.}},
\bauthor{\bsnm{Prelic},~\bfnm{A.}\binits{A.}},
\bauthor{\bsnm{Von~Rohr},~\bfnm{P.}\binits{P.}},
\bauthor{\bsnm{Thiele},~\bfnm{L.}\binits{L.}} \betal{et~al.}
(\byear{2004}).
\btitle{{Sparse graphical Gaussian modeling of the isoprenoid gene network in
Arabidopsis thaliana}}.
\bjournal{Genome Biology}
\bvolume{5}
\bpages{1--13}.
\bptok{imsref}%
\end{barticle}
\endbibitem

\bibitem[\protect\citeauthoryear{Xue and Zou}{2011a}]{xue2011}
\begin{bmisc}[author]
\bauthor{\bsnm{Xue},~\bfnm{L.}\binits{L.}} \AND
\bauthor{\bsnm{Zou},~\bfnm{H.}\binits{H.}}
(\byear{2011a}).
\bhowpublished{On estimating sparse correlation matrices of semiparametric Gaussian
copulas. Technical report, Univ. Minnesota.}
\bptok{imsref}%
\end{bmisc}
\endbibitem



\bibitem[\protect\citeauthoryear{Xue and Zou}{2011b}]{XZ11b}
%
\begin{bmisc}[auto]
\bauthor{\bsnm{Xue},~\bfnm{L.}\binits{L.}} \AND
\bauthor{\bsnm{Zou},~\bfnm{H.}\binits{H.}}
(\byear{2011b}).
\bhowpublished{Regularized rank-based estimation of high-dimensional nonparanormal graphical models.
Technical report, Univ. Minnesota}.
\bptok{imsref}%
\end{bmisc}
%
\endbibitem


\bibitem[\protect\citeauthoryear{Xue and Zou}{2012}]{supp}
\begin{bmisc}[auto]
\bauthor{\bsnm{Xue},~\bfnm{L.}\binits{L.}} \AND
\bauthor{\bsnm{Zou},~\bfnm{H.}\binits{H.}}
(\byear{2012}).
\bhowpublished{Supplement to ``Regularized rank-based estimation of high-dimensional
nonparanormal graphical models.'' DOI:\doiurl{10.1214/12-AOS1041SUPP}.}
\bptok{imsref}%
\end{bmisc}
\endbibitem

\bibitem[\protect\citeauthoryear{Yuan}{2010}]{yuan2010}
\begin{barticle}[mr]
\bauthor{\bsnm{Yuan},~\bfnm{Ming}\binits{M.}}
(\byear{2010}).
\btitle{High dimensional inverse covariance matrix estimation via linear
programming}.
\bjournal{J. Mach. Learn. Res.}
\bvolume{11}
\bpages{2261--2286}.
\bid{issn={1532-4435}, mr={2719856}}
\bptok{imsref}%
\end{barticle}
\endbibitem

\bibitem[\protect\citeauthoryear{Yuan and Lin}{2007}]{yuan2007}
\begin{barticle}[mr]
\bauthor{\bsnm{Yuan},~\bfnm{Ming}\binits{M.}} \AND
\bauthor{\bsnm{Lin},~\bfnm{Yi}\binits{Y.}}
(\byear{2007}).
\btitle{Model selection and estimation in the {G}aussian graphical model}.
\bjournal{Biometrika}
\bvolume{94}
\bpages{19--35}.
\bid{doi={10.1093/biomet/asm018}, issn={0006-3444}, mr={2367824}}
\bptok{imsref}%
\end{barticle}
\endbibitem

\bibitem[\protect\citeauthoryear{Zhou et~al.}{2011}]{zhou2010}
\begin{barticle}[mr]
\bauthor{\bsnm{Zhou},~\bfnm{Shuheng}\binits{S.}},
\bauthor{\bsnm{R{\"u}timann},~\bfnm{Philipp}\binits{P.}},
\bauthor{\bsnm{Xu},~\bfnm{Min}\binits{M.}} \AND
\bauthor{\bsnm{B{\"u}hlmann},~\bfnm{Peter}\binits{P.}}
(\byear{2011}).
\btitle{High-dimensional covariance estimation based on {G}aussian graphical
models}.
\bjournal{J. Mach. Learn. Res.}
\bvolume{12}
\bpages{2975--3026}.
\bid{issn={1532-4435}, mr={2854354}}
\bptok{imsref}%
\end{barticle}
\endbibitem

\bibitem[\protect\citeauthoryear{Zou}{2006}]{zou2006}
\begin{barticle}[mr]
\bauthor{\bsnm{Zou},~\bfnm{Hui}\binits{H.}}
(\byear{2006}).
\btitle{The adaptive lasso and its oracle properties}.
\bjournal{J. Amer. Statist. Assoc.}
\bvolume{101}
\bpages{1418--1429}.
\bid{doi={10.1198/016214506000000735}, issn={0162-1459}, mr={2279469}}
\bptok{imsref}%
\end{barticle}
\endbibitem

\end{thebibliography}
\end{document}